\documentclass[12pt]{amsart}
\usepackage[all,dvips,arc,curve,color,frame]{xy}
\usepackage{graphics}
\usepackage{amsmath,amssymb}
\usepackage{color}
\usepackage{url}
\usepackage{epsfig}
\title{$4$-regular and self-dual analogs of fullerenes}

\author{Mathieu Dutour Sikiri\'c}
\address{M. Dutour Sikiri\'c, Rudjer Boskovi\'c Institute, Bijenicka 54, 10000 Zagreb, Croatia}
\email{mdsikir@irb.hr}

\author{Michel Deza}
\address{M. Deza, JAIST, Nomi, Ishikawa-ken, 923-1292, Japan}
\email{Michel.Deza@ens.fr}

\def\QuotS#1#2{\leavevmode\kern-.0em\raise.2ex\hbox{$#1$}\kern-.1em/\kern-.1em\lower.25ex\hbox{$#2$}}

\begin{document}
\newcommand{\KK}{\ensuremath{\mathbb{K}}}
\newcommand{\RR}{\ensuremath{\mathbb{R}}}
\newcommand{\NN}{\ensuremath{\mathbb{N}}}
\newcommand{\QQ}{\ensuremath{\mathbb{Q}}}
\newcommand{\ZZ}{\ensuremath{\mathbb{Z}}}

\newtheorem{theorem}{Theorem}
\newtheorem{proposition}[theorem]{Theorem}
\newtheorem{corollary}[theorem]{Corollary}
\newtheorem{lemma}[theorem]{Lemma}
\newtheorem{problem}[theorem]{Problem}
\newtheorem{conjecture}[theorem]{Conjecture}
\newtheorem{claim}[theorem]{Claim}
\newtheorem{remark}[theorem]{Remark}
\newtheorem{definition}[theorem]{Definition}

\begin{abstract}
An {\em $i$-hedrite} 
is a $4$-regular plane graph with faces of size $2$, 
$3$ and $4$. 
We do a short survey of their known properties \cite{oct2,DHL,octa} and
explain some new algorithms that allow their efficient enumeration.
Using this we give the symmetry groups of all $i$-hedrites and the minimal
representative for each.
We also review the link of $4$-hedrites with knot theory and the classification
of $4$-hedrites with simple central circuits.

An {\em $i$-self-hedrite} is a self-dual
plane graph with faces and vertices of size/degree $2$, $3$ and $4$.
We give a new efficient algorithm for enumerating them based on $i$-hedrites.
We give a classification of their possible symmetry groups and a
classification of $4$-self-hedrites of symmetry $T$, $T_d$ in terms of 
the {\em Goldberg-Coxeter construction}.
Then we give a method for enumerating $4$-self-hedrites with simple
zigzags.
\end{abstract}

\maketitle

\section{Introduction}
A {\em fullerene} is a $3$-regular  plane graph whose faces have size $5$ or $6$.
As a consequence of Euler's formula any fullerene has exactly $12$ $5$-gonal faces.
For a $3$-regular  plane graph $G$ and a $r$-gonal face $F$ of $G$, the quantity
$6-r$ is called {\em curvature} and Euler's formula is then
a statement about the curvature on the sphere.
A natural generalization of fullerene is the class of $3$-regular  
plane graphs with faces of size between $3$ and $6$ (see, for example,
\cite{Group345}).

Here we consider another generalization, that is  a suitable $k$-regular plane graph. The Euler formula $V-E+F=2$ becomes then 
\begin{equation*}
\sum_{j=2}^{\infty}p_j(s - j)=\frac{4k}{k-2}\mbox{~with~}s=\frac{2k}{k-2};
\end{equation*}
we will permit $2$-gons (doubled edges) but not $1$-gons.
The only integral pairs $(s,k)$ are $(6,3)$, $(4,4)$ and $(3,6)$.
We will permit only $s$- and $(s-1)$-gonal faces.
So, $p_{s-1}=\frac{4k}{k-2}$ and $p_s$ is not bounded.
The number $n$ of vertices is 
\begin{equation*}
n=\frac{4(k+2)}{{(k-2)^2}}+p_s\frac{2}{k-2}.
\end{equation*}
For $k=3,4$ and $6$ we get spherical analogs of the
regular partition of the Euclidean plane $E^2$: $\{6^3\}$, $\{4^4\}$ and $\{3^6\}$, 
respectively,
where $12$ pentagons, $8$
triangles and $6$ doubled edges play role of ``defects'', disclinations needed to increase the
curvature to the one of sphere $S^2$.
The graphs with smallest number $n$ of vertices have only
$(s-1)$-gons; they are Dodecahedron, Octahedron and $Bundle_6$
($2$ vertices connected by $6$ edges)
for  $k=3,4,6$, respectively.
The case $k=3$ gives fullerenes.
The case $k=4$, i.e., of $4$-regular plane graphs with faces of
size $3$ or $4$, gives {\em octahedrites} treated in the foundational
paper \cite{octa}.
Let us call graphs in the remaining case
$k=6$ ($6$-regular plane graphs with faces of size $2$ or $3$) {\em bundelites}.
Thurston's work \cite{T} implies that fullerenes can be parametrized by $10$ Eisenstein integers and the number of  fullerenes with $n$ vertices
grows as $n^9$; those results can be generalized to octahedrites and bundelites.
The ring of definition for octahedrites, respectively bundelites, is the
Gaussian, respectively Eisenstein integers. Those cases belong to two of $94$
cases enumerated in \cite{T}.

We present here a short review of known facts about octahedrites
as established in \cite{oct2,DHL,octa} and present a few new facts and
applications.
We give the possible symmetry groups of octahedrites and the graphs of minimal
vertex-sets realizing them.
Then we show how octahedrites can be used for the enumeration of a 
more general {\em $i$-hedrites}, i.e. $4$-regular plane 
graphs with faces of size $2$, $3$ or $4$ and $p_2 + p_3=i$.

Then we consider {\em central circuit partition} of the edge-sets of
octahedrites and the corresponding knot-theoretic notions, that is
{\em alternating knot}, {\em Borromean link}, and {\em equivalence}.

A {\em $i$-self-hedrite} is a plane graph
with vertices and faces of size
$2$, $3$ or $4$ that is isomorphic to its dual with $p_2+p_3=i$.
Such graphs have $2p_2 + p_3=4$ and $i$-self-hedrites can be enumerated
effectively by using $2i$-hedrites with a method detailed below.
We determine their possible symmetry groups and we list the minimal
representatives for each of them.
We characterize the $4$-self-hedrites of symmetry $T$ or $T_d$ in terms
of the {\em Goldberg-Coxeter construction} for octahedrites.
Then we give a method based on $2i$-hedrites for determining the
$i$-self-hedrites with {\em simple zigzags}.

The computations of this paper were done using the GAP computer algebra system
and the computer packages {\tt polyhedral}, {\tt plangraph} of the first
author. 
The enumeration of octahedrites was done using the {\tt ENU} program by 
O. Heidemeier \cite{He,Heidemeier} and the program {\tt CaGe} \cite{CaGe} was used
for making the drawings.

\section{Structural properties}

A {\em plane graph} is a graph drawn on the plane with edges intersecting only at
vertices. A graph $G$ is {\em $3$-connected} if after removing any $2$ vertices
of $G$ the resulting graph is connected. A {\em $3$-polyhedron} is a $3$-dimensional
polytope, its skeleton defines a $3$-connected plane graph and it is known
that this characterizes the skeleton of $3$-polytopes.
Furthermore \cite{Mani71}, a $3$-connected plane graph $G$ can be
represented as a skeleton of a $3$-polytope $P$ such that any symmetry
of $G$ is realized as a isometry of the polytope $P$.
We refer to \cite{book3} for more details on such questions.

It is proved in \cite{oct2} that any octahedrite is $3$-connected
which implies that its symmetry groups is realized as isometry of $3$-space.
Since those group have been classified long ago and are much used in
chemistry, we can use the chemical nomenclature here (see, for a possible
presentation, \cite{pointgroup}).

An octahedrite exists for any $n\geq 6$ except $n=7$ (see \cite[page 282]{Gr1}).
For a $4$-regular graph with $p_j$ denoting the number of faces of size
$j$, the classical Euler formula $V-E+F=2$ can be rewritten (see \cite[Chapter 1]{book3}, for the easy details) as
\begin{equation}\label{EulerFormula4valent}
\sum_{j=2}^{\infty} (4-j) p_j=8.
\end{equation}
For octahedrites this directly implies $p_3=8$.
Octahedron is the unique octahedrite with $n=6$.

\begin{figure}
\begin{center}
\begin{minipage}[b]{2.3cm}
\centering
\epsfig{height=16mm, file=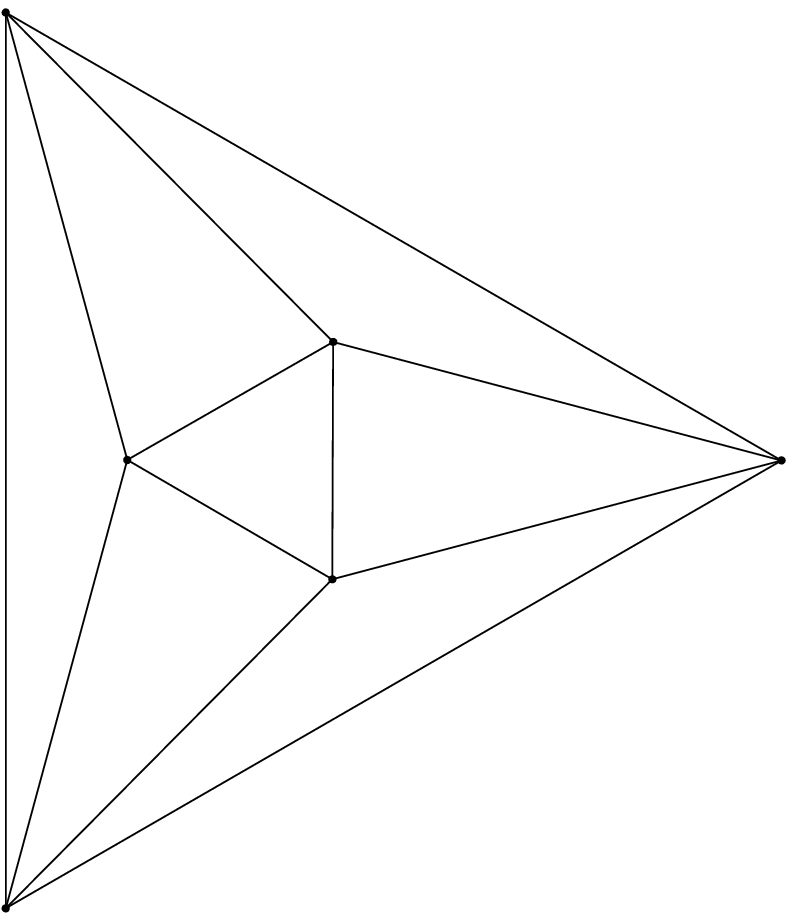}\par
$O_h$, $n=6$
\end{minipage}
\begin{minipage}[b]{2.3cm}
\centering
\epsfig{height=16mm, file=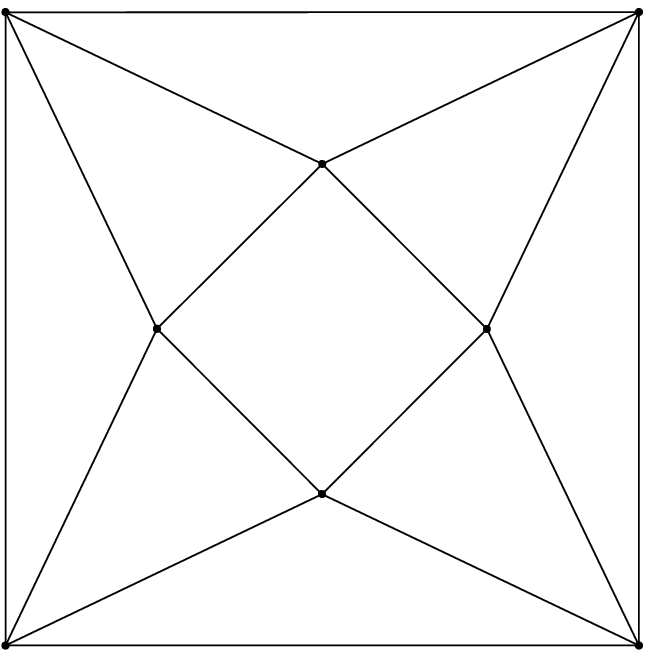}\par
$D_{4d}$, $n=8$
\end{minipage}
\begin{minipage}[b]{2.3cm}
\centering
\epsfig{height=16mm, file=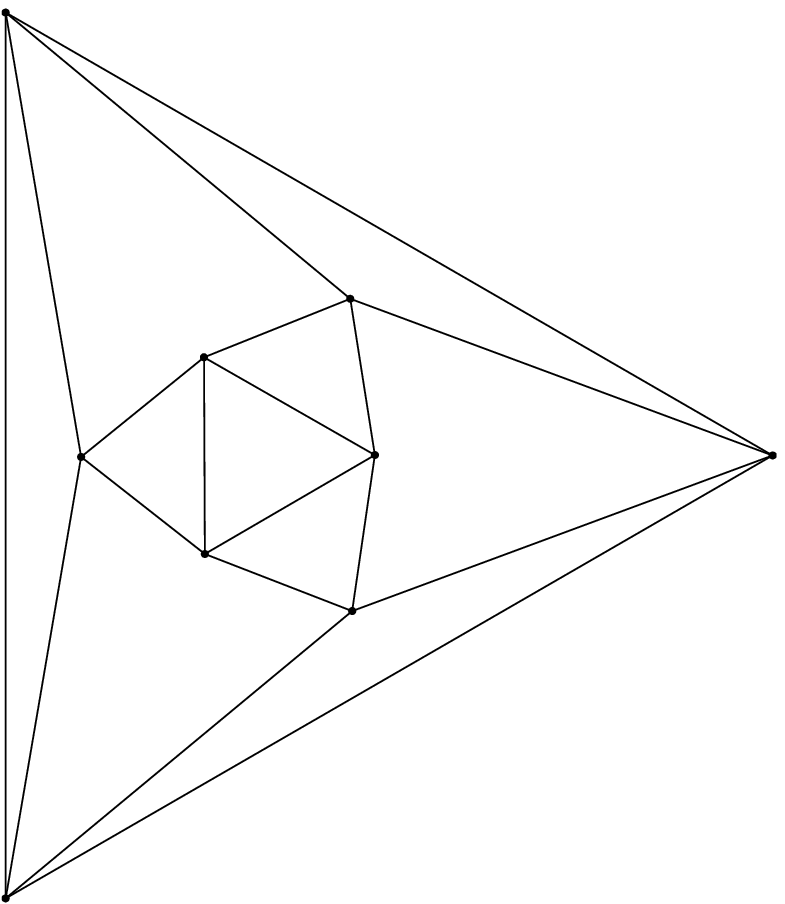}\par
$D_{3h}$, $n=9$
\end{minipage}
\begin{minipage}[b]{2.3cm}
\centering
\epsfig{height=16mm, file=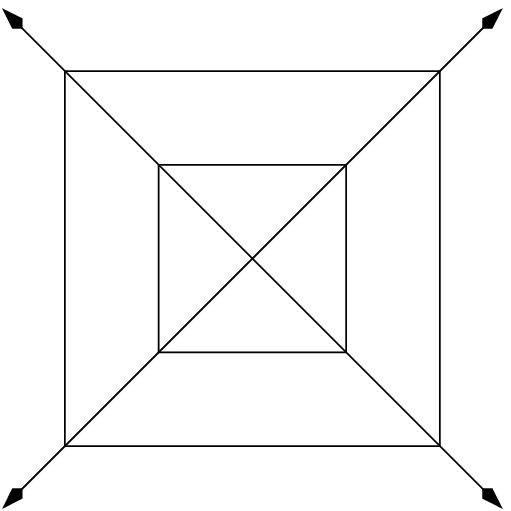}\par
$D_{4h}$, $n=10$
\end{minipage}
\begin{minipage}[b]{2.3cm}
\centering
\epsfig{height=16mm, file=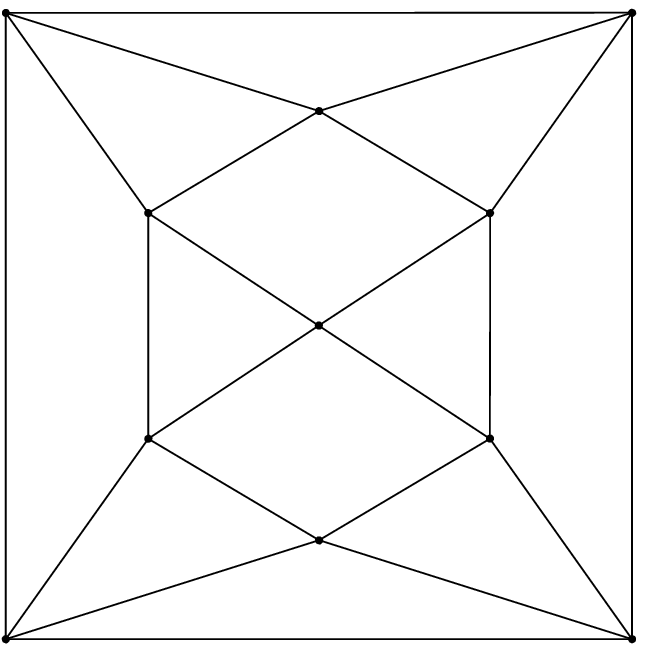}\par
$C_{2v}$, $n=11$
\end{minipage}
\begin{minipage}[b]{2.3cm}
\centering
\epsfig{height=16mm, file=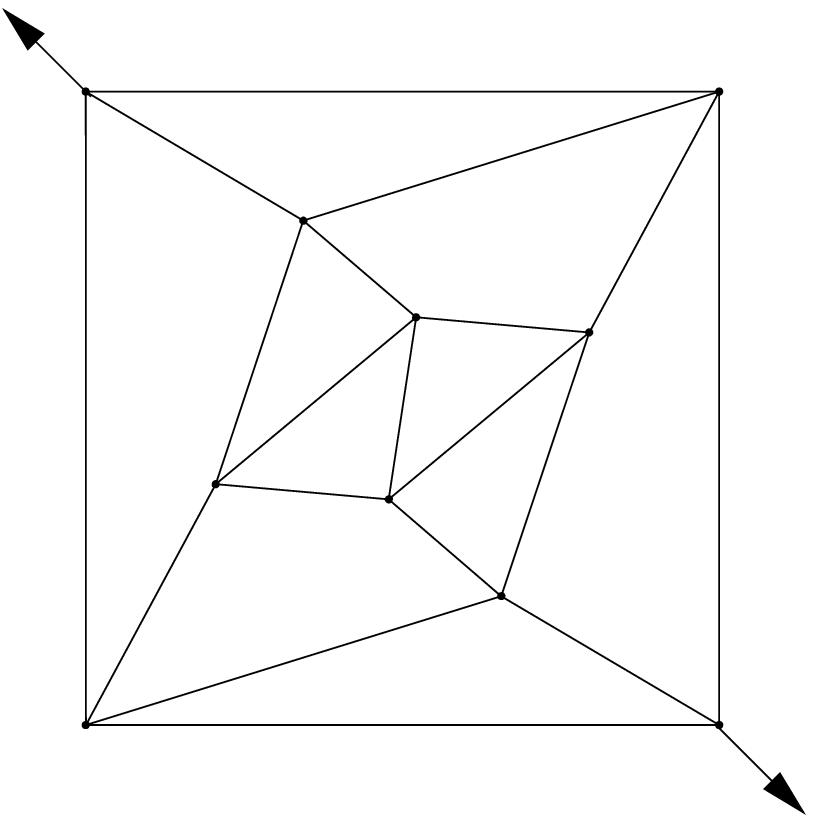}\par
$D_{2}$, $n=12$
\end{minipage}
\begin{minipage}[b]{2.3cm}
\centering
\epsfig{height=16mm, file=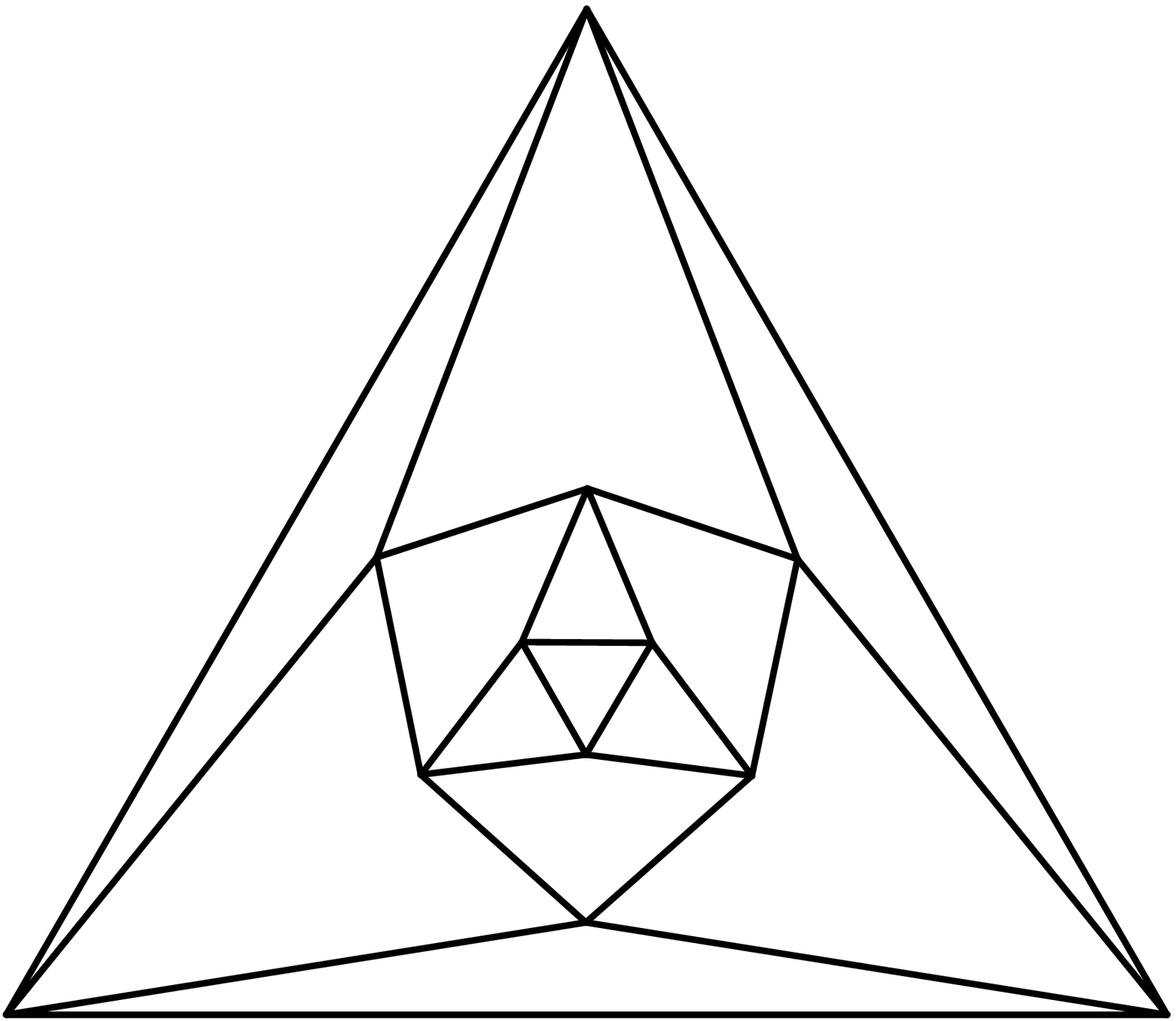}\par
$D_{3d}$, $n=12$
\end{minipage}
\begin{minipage}[b]{2.3cm}
\centering
\epsfig{height=16mm, file=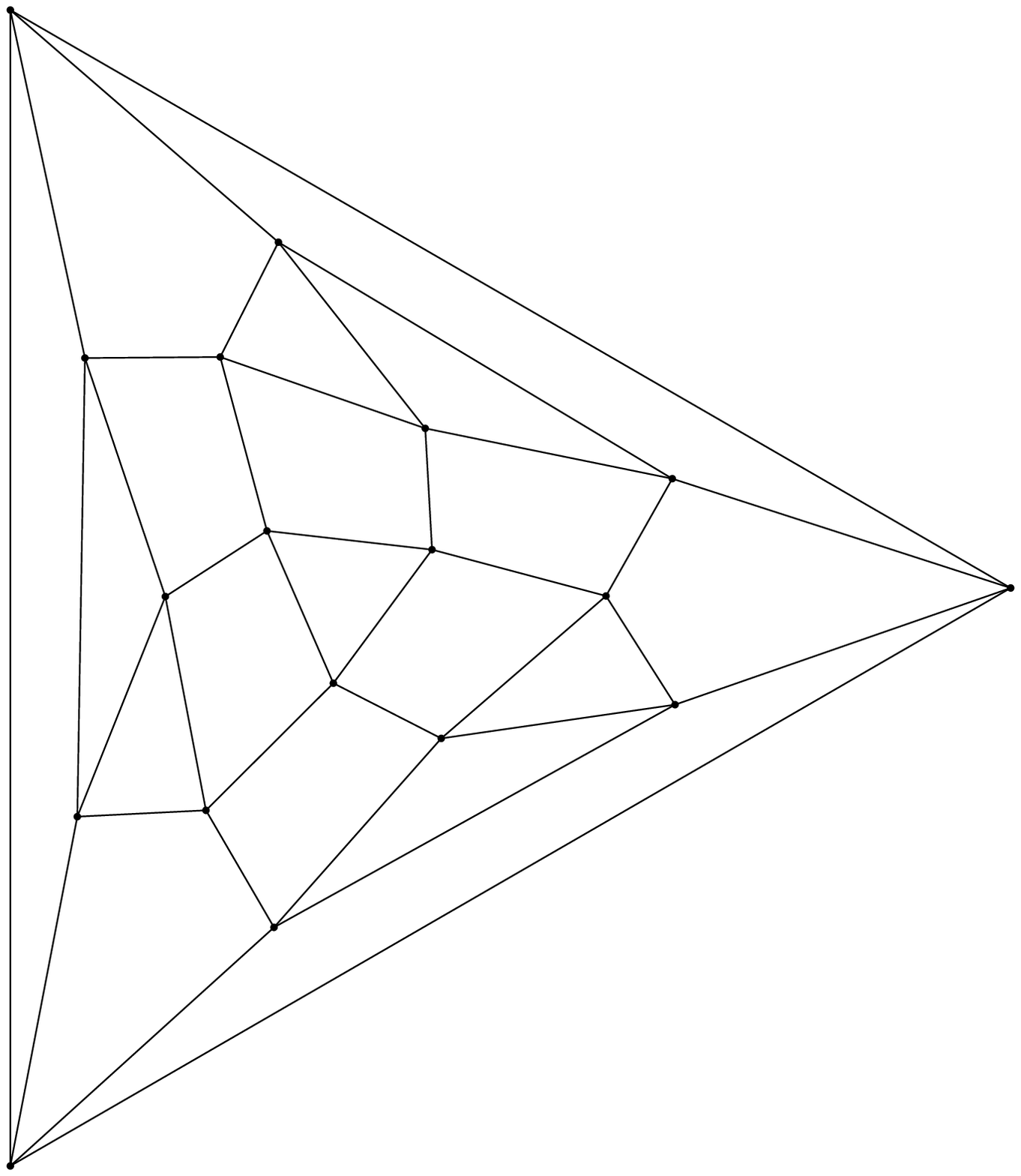}\par
$D_{3}$, $n=12$
\end{minipage}
\begin{minipage}[b]{2.3cm}
\centering
\epsfig{height=16mm, file=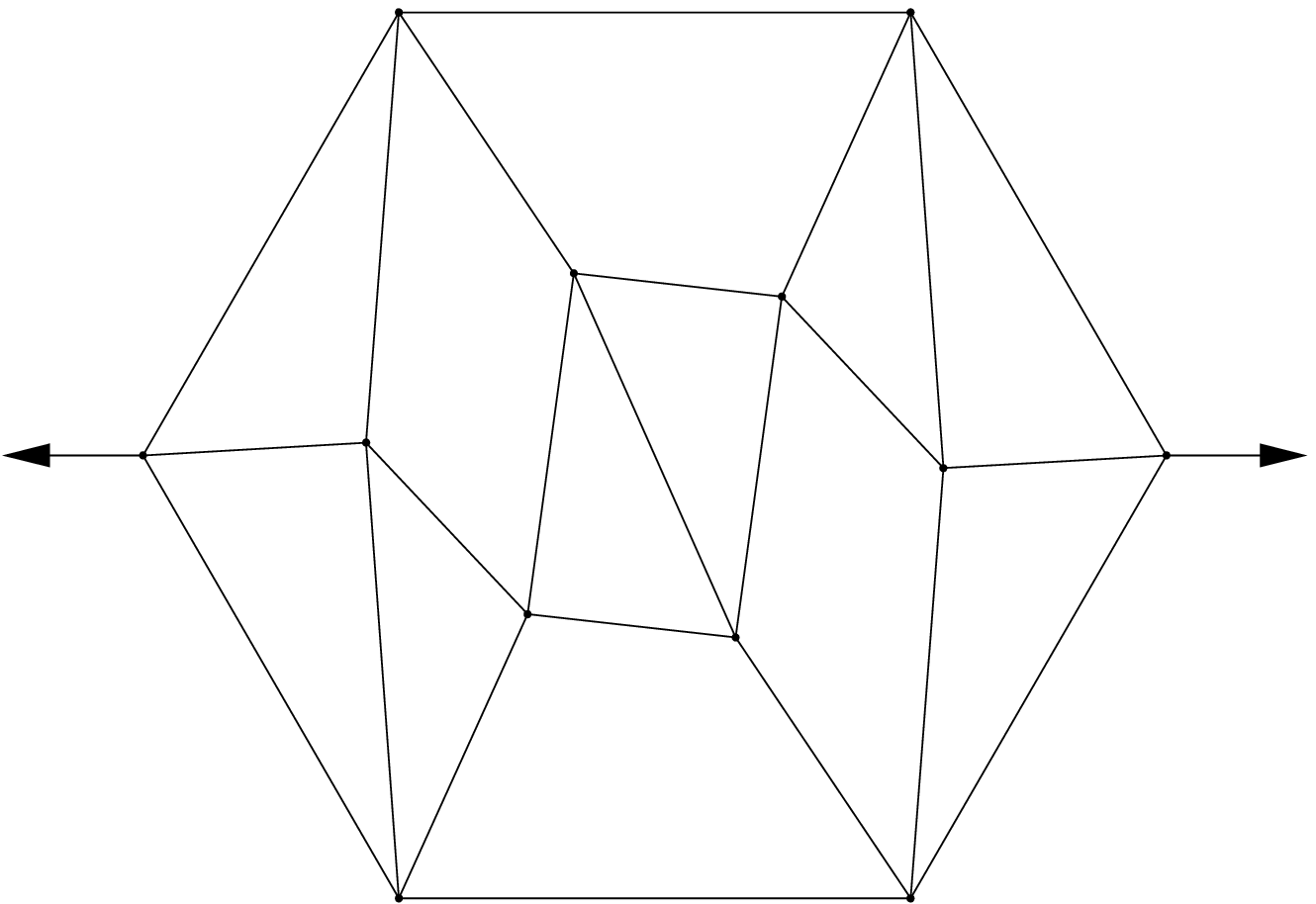}\par
$C_{2}$, $n=12$
\end{minipage}
\begin{minipage}[b]{2.3cm}
\centering
\epsfig{height=16mm, file=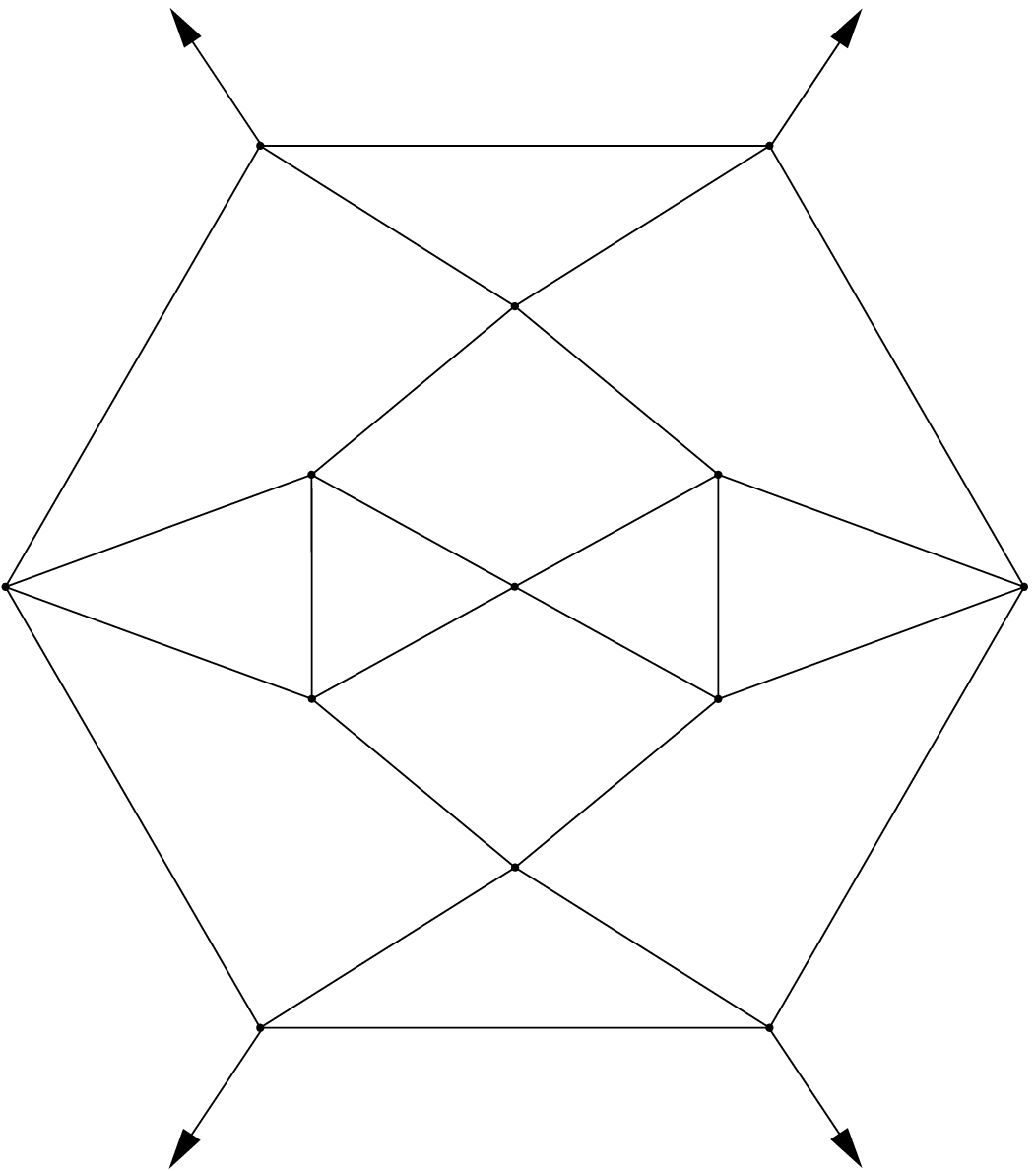}\par
$D_{2d}$, $n=14$
\end{minipage}
\begin{minipage}[b]{2.3cm}
\centering
\epsfig{height=16mm, file=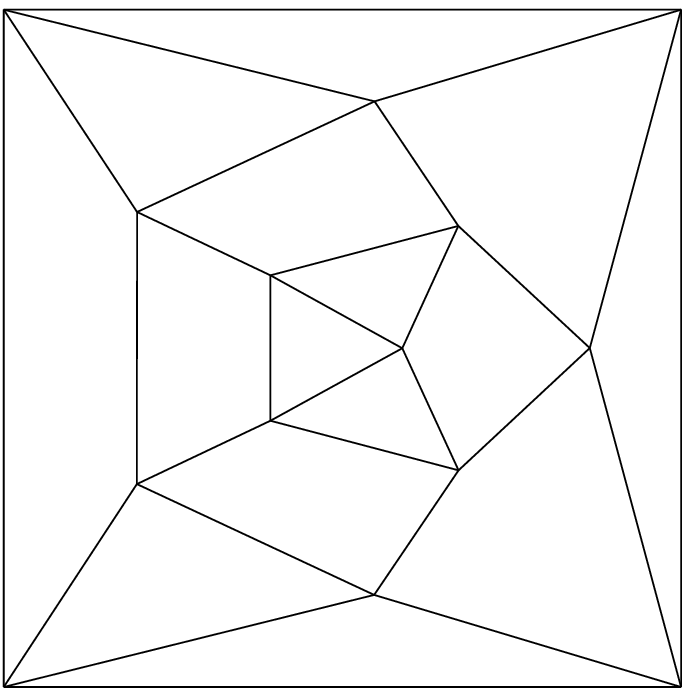}\par
$C_{s}$, $n=14$
\end{minipage}
\begin{minipage}[b]{2.3cm}
\centering
\epsfig{height=16mm, file=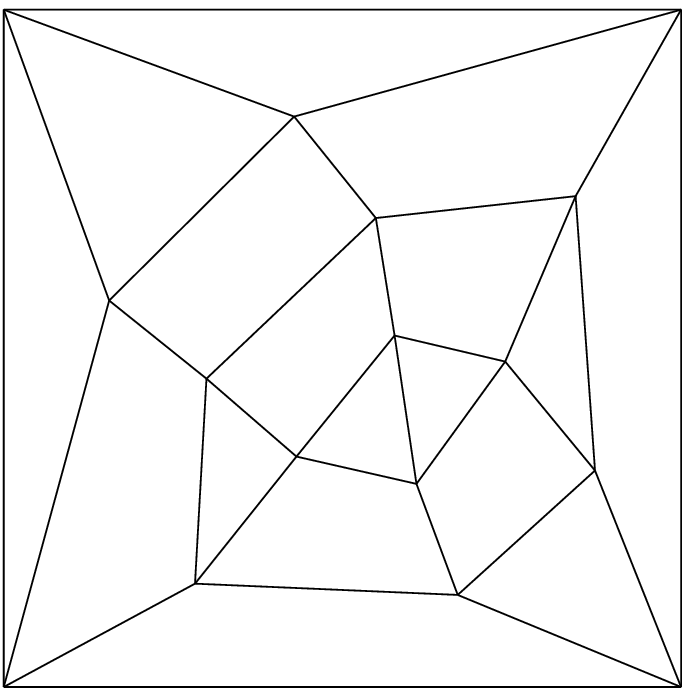}\par
$C_{1}$, $n=16$
\end{minipage}
\begin{minipage}[b]{2.3cm}
\centering
\epsfig{height=16mm, file=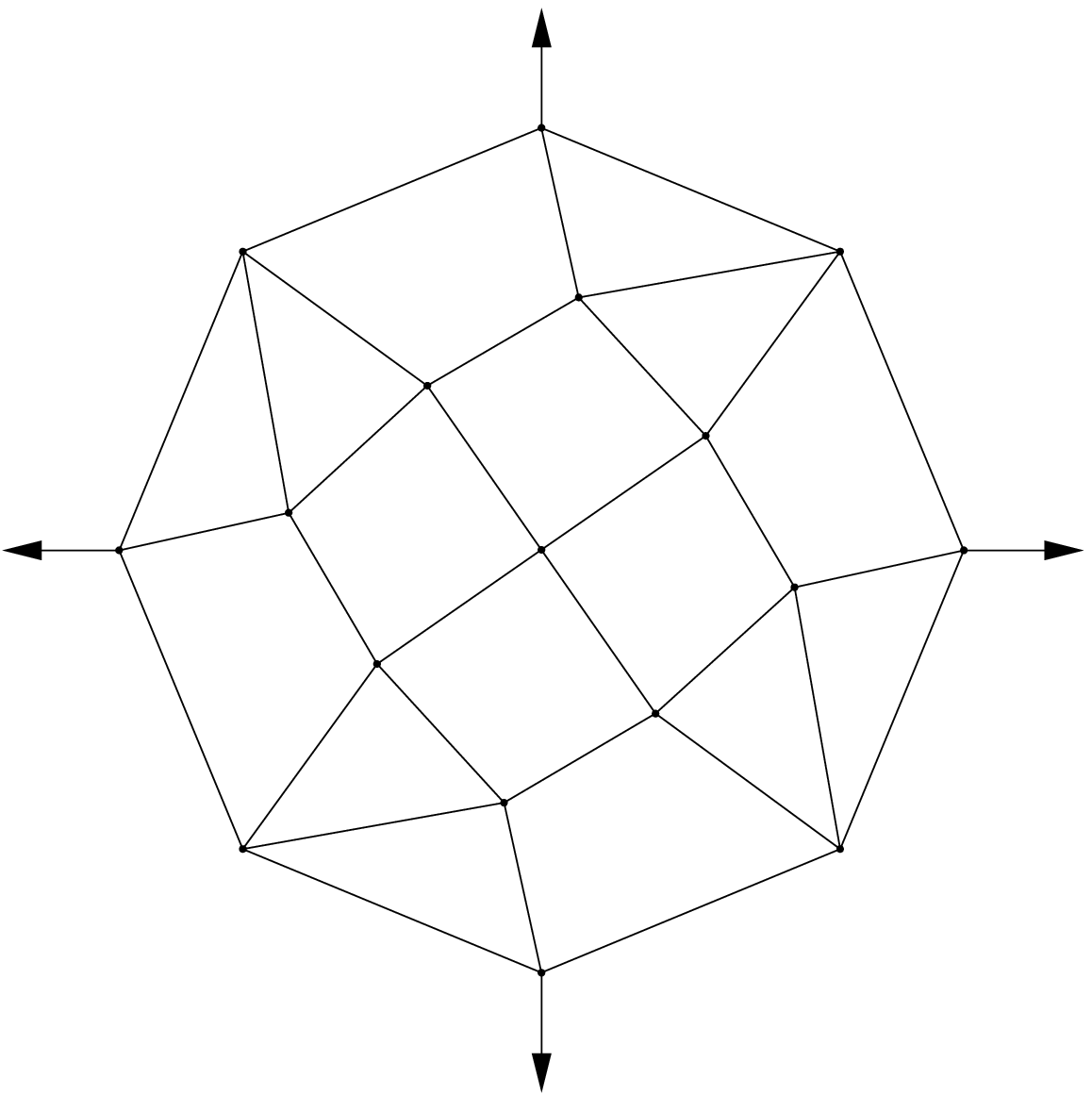}\par
$D_{4}$, $n=18$
\end{minipage}
\begin{minipage}[b]{2.3cm}
\centering
\epsfig{height=16mm, file=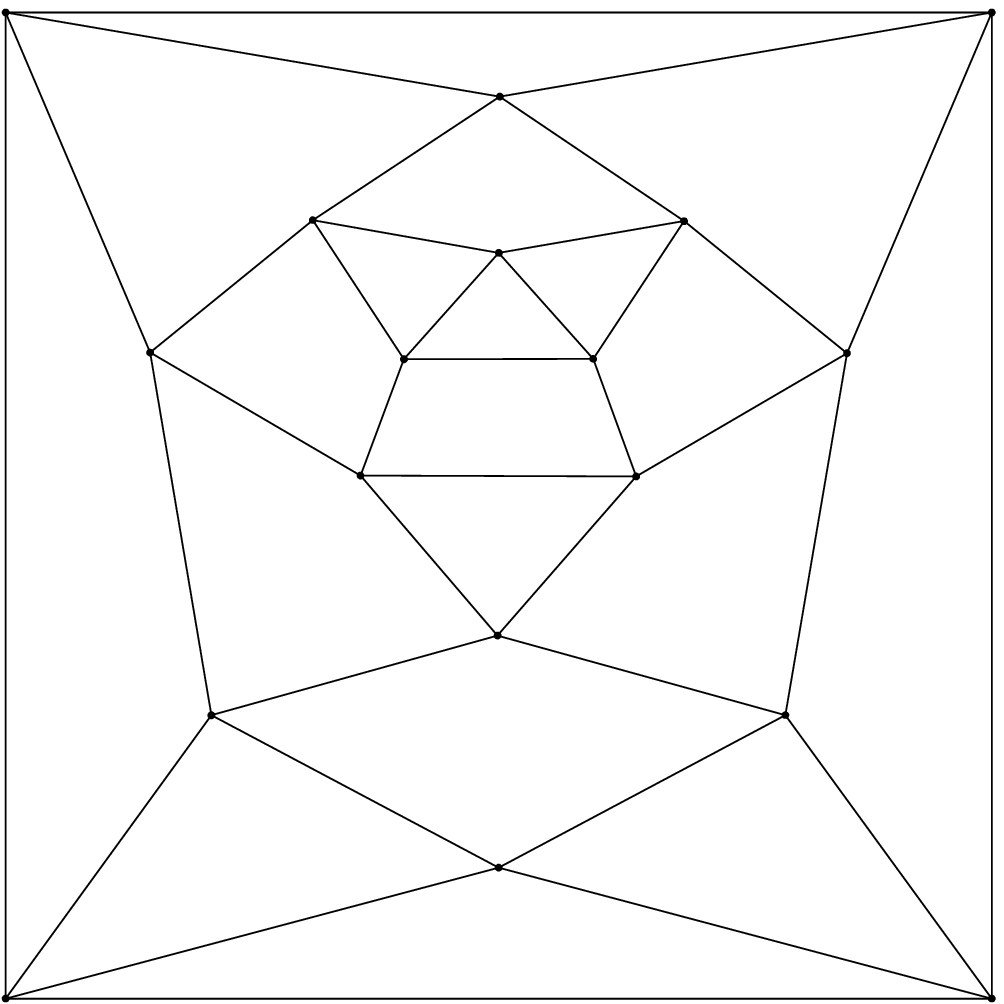}\par
$C_{2h}$, $n=18$
\end{minipage}
\begin{minipage}[b]{2.3cm}
\centering
\epsfig{height=16mm, file=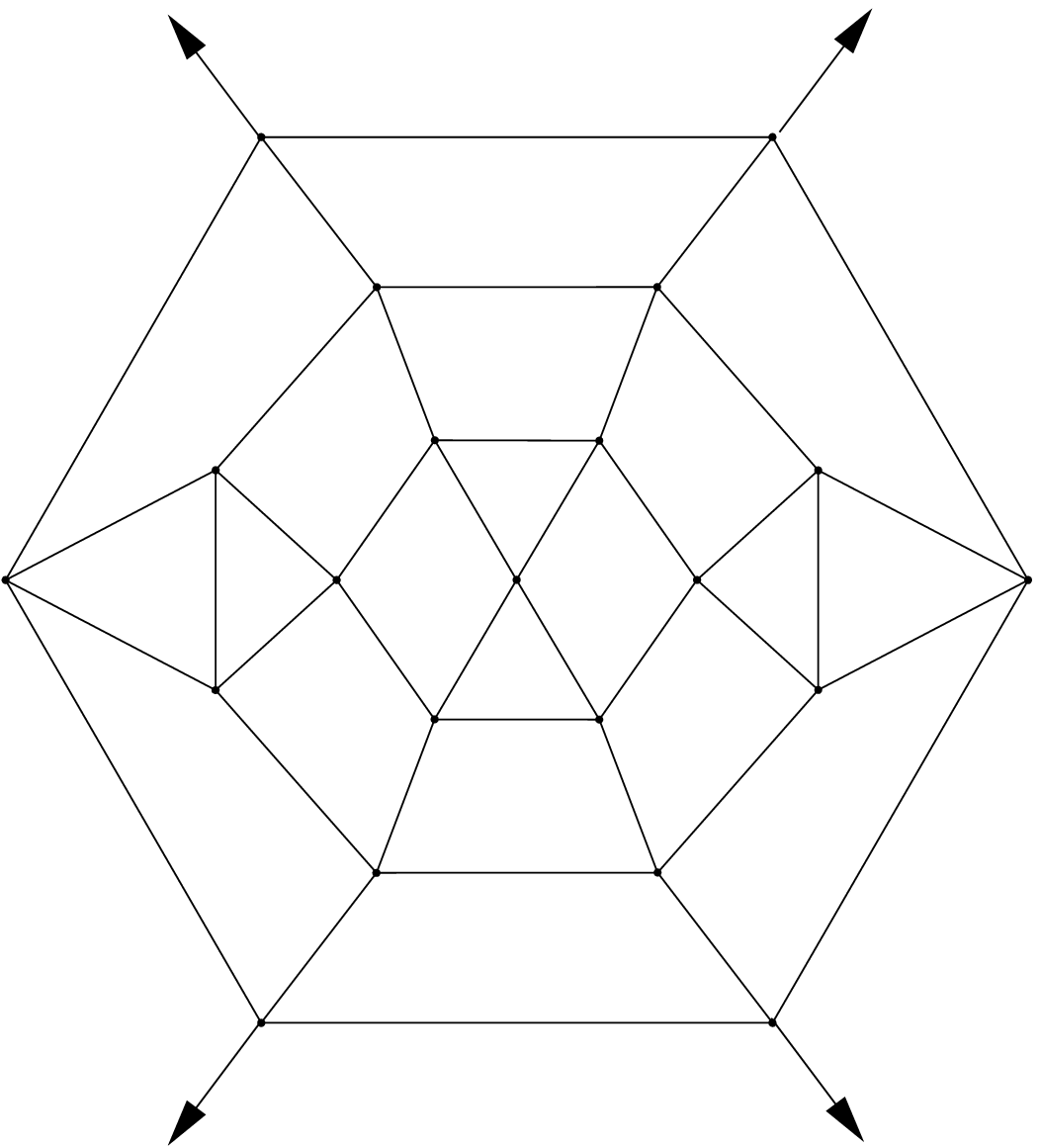}\par
$D_{2h}$, $n=22$
\end{minipage}
\begin{minipage}[b]{2.3cm}
\centering
\epsfig{height=16mm, file=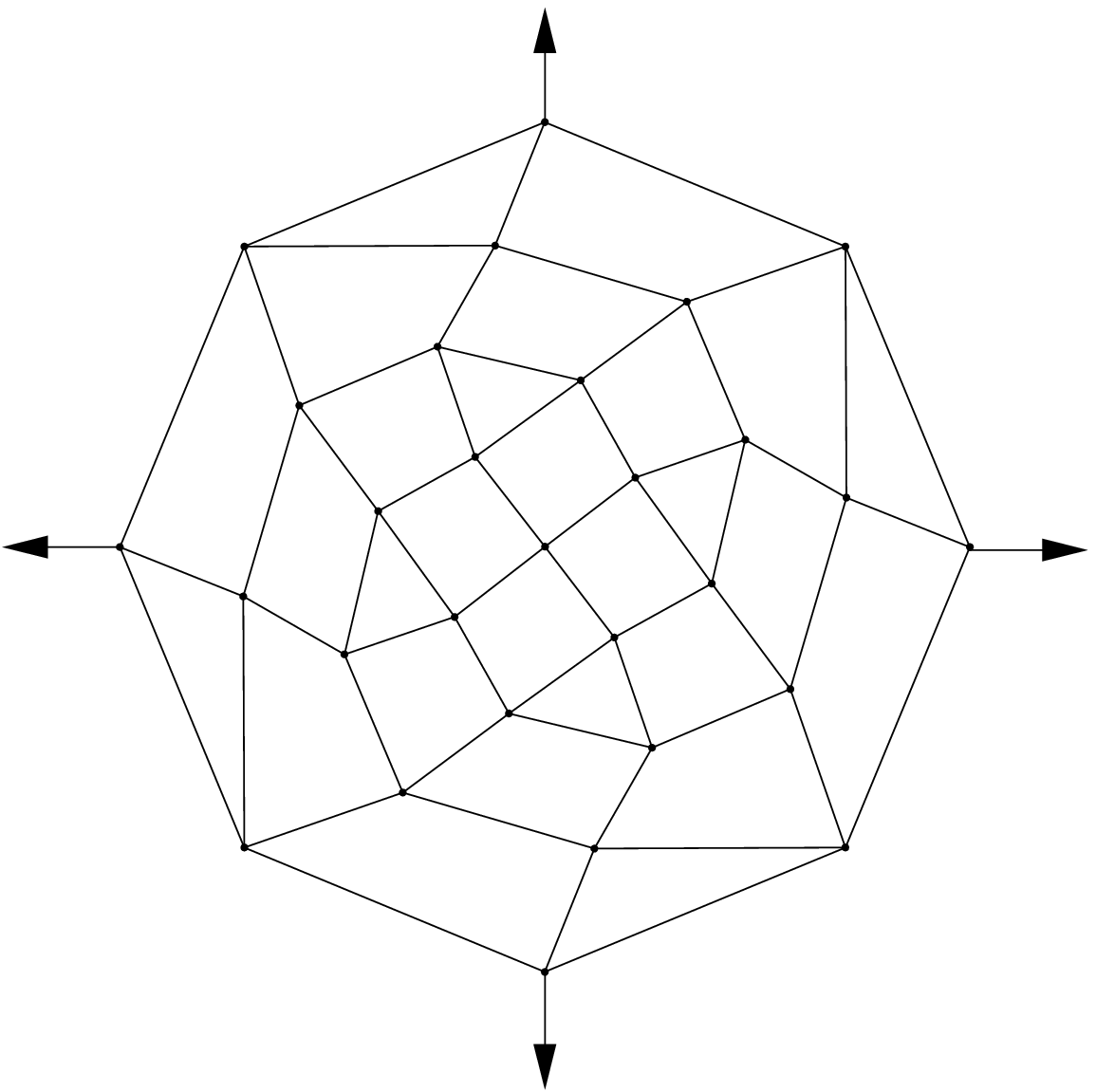}\par
$O$, $n=30$
\end{minipage}
\begin{minipage}[b]{2.3cm}
\centering
\epsfig{height=16mm, file=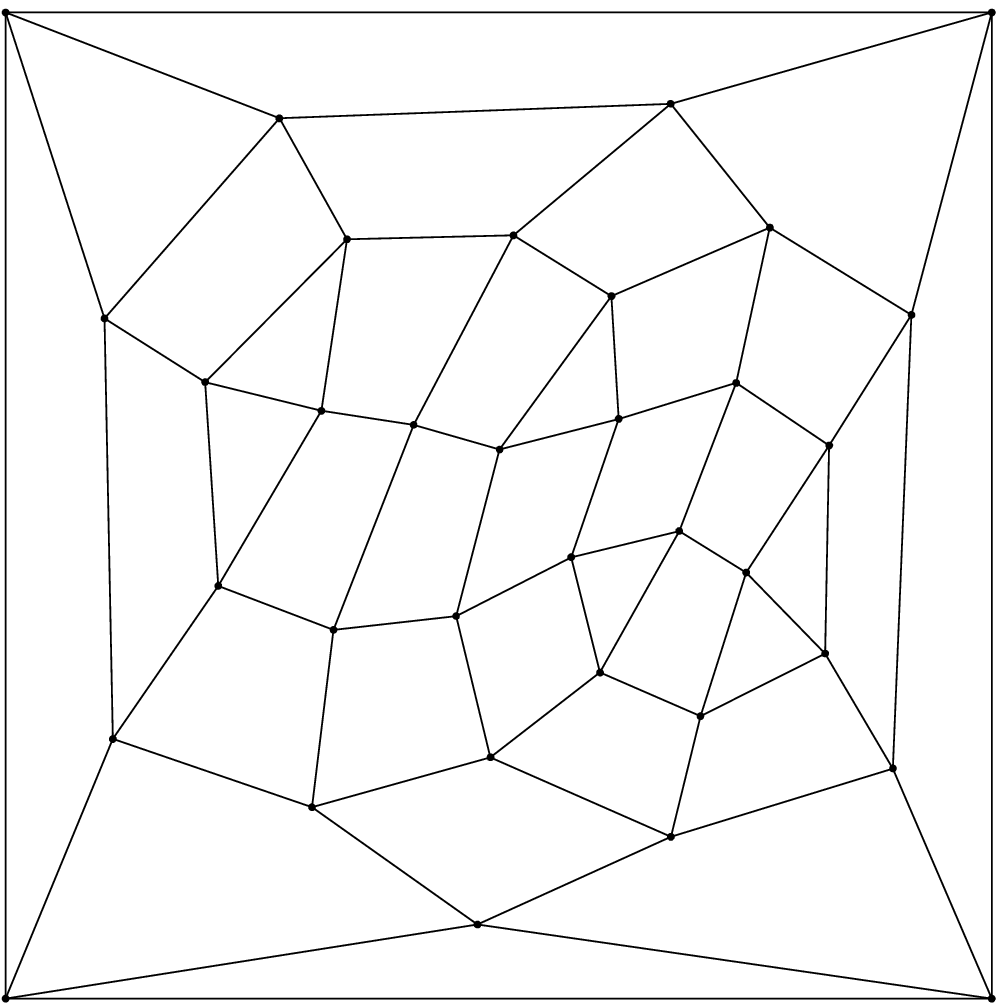}\par
$C_{i}$, $n=34$
\end{minipage}
\begin{minipage}[b]{2.3cm}
\centering
\epsfig{height=16mm, file=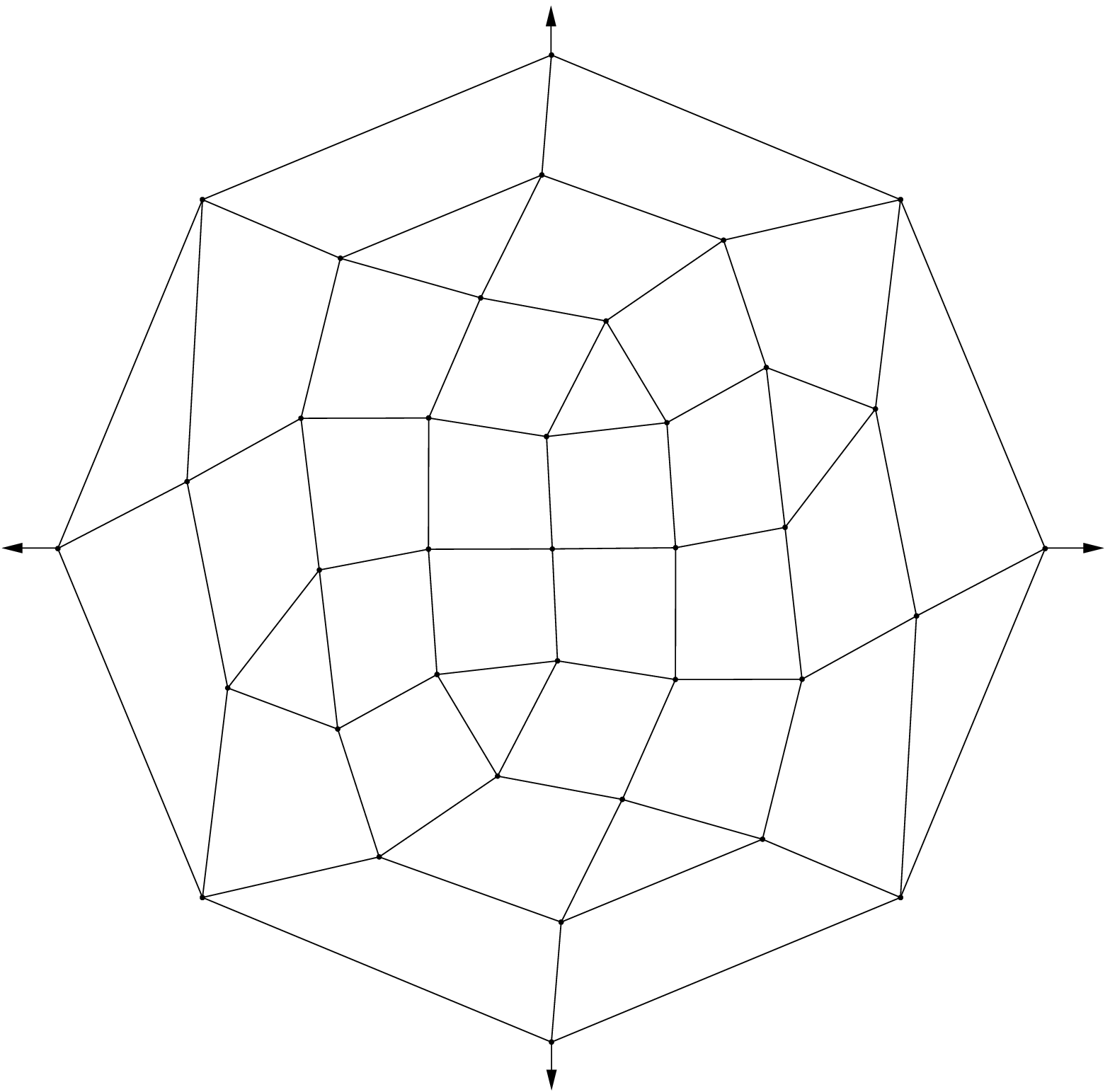}\par
$S_{4}$, $n=38$
\end{minipage}

\end{center}
\caption{Minimal representatives for each possible symmetry group of an octahedrite}
\label{MinimalRepresentative}
\end{figure}

\begin{theorem}\label{GroupsOfOctahedrites}
The only symmetry groups of octahedrites are: $C_{1}$, $C_s$, $C_2$, $C_{2v}$, $C_i$, $C_{2h}$, $S_4$, $D_2$, $D_{2d}$, $D_{2h}$, $D_3$, $D_{3d}$, $D_{3h}$, $D_4$, $D_{4d}$, $D_{4h}$, $O$, $O_h$.
The minimal possible representative are given in Figure \ref{MinimalRepresentative}.
\end{theorem}
The proof that the list of groups is complete is given in \cite{oct2}, but the
minimal possible representatives were not determined at the time.
The method is first to go through the restrictions that vertex degree and face size impose.
An $m$-fold axis of rotation has necessarily $m=4$ (passing though a face of size $4$ or a vertex), $m=3$ (axis passing through a face of size $3$), or $m=2$ (axis passing though an edge, a vertex of degree $4$ or a face of size $4$).
Then the classification of point groups gives a list of possible candidates.
Some candidates are excluded for reason of orbit size and other similar simple
arguments. But some groups
are excluded for a subtler reason: the existence of a symmetry implies another
symmetry. For example a $3$-, $4$-fold axis of symmetry, i.e. $C_3$, $C_4$
implies actually at least $D_3$, $D_4$ for possible symmetry groups.
See \cite{oct2} for details.

On the other hand, finding the minimal possible representative is done in a very non-clever way: we look at all the generated octahedrites and select the representatives with minimal vertex-sets.
The enumeration of octahedrites was done by using the program {\tt ENU} (see \cite{He,Heidemeier}) by O. Heidemeier, that enumerates classes of $4$-regular graphs with constraint on the size of their faces, fairly efficiently.

\section{Generation of $i$-hedrites}

Define an {\em $i$-hedrite} to be a $4$-regular $n$-vertex plane graph, whose 
faces have size $2$, $3$ and $4$ only and $p_2+p_3=i$
(see, for more details, \cite{oct2}).
Using Formula \eqref{EulerFormula4valent}, we get for an $i$-hedrite
$2p_2 + p_3=8$ and the only solutions 
are $i=4,5,6,7,8$, which have, respectively, $(p_2, p_3)$=$(4,0)$, $(3,2)$, $(2,4)$, $(1,6)$ and $(0,8)$.
So, $8$-hedrites are octahedrites.
We will be concerned here only about the generation of $i$-hedrites.
Actually, $4$-hedrites admit a reasonably simple explicit description,
see \cite{octa} and \cite[Chapter 2]{book3}.
So, it remains to find efficient methods for the enumeration of $5$-, 
$6$- and $7$-hedrites. The program {\tt ENU} cannot deal with faces
of size $2$; so, we sought a method that allows for reasonable enumeration
of such graphs. See Table \ref{EnumerationIHedrite} for
the number of $i$-hedrites with at most $70$ vertices.
Easy to check that an $n$-vertex $i$-hedrite exists for even $n\ge 2$ if
$i=4$, $n\ge 5$ (and $n=3$) if $i=5$,  $n\ge 4$ if $i=6$, $n\ge 7$ if
$i=7$,  $n\ge 8$ (and $n=6$) if $i=8$.

\begin{table}
\caption{Number of $i$-hedrites, $4\le i\le 8$,  with  $2\leq n\leq 70$}
\label{EnumerationIHedrite}
{\scriptsize
\begin{equation*}
\begin{array}{||c|c|c|c|c|c||c|c|c|c|c|c||c|c|c|c|c|c||}
\hline
\hline
n & {\bf 4} & {\bf 5} & {\bf 6} & {\bf 7} & {\bf 8} & n & {\bf 4} & {\bf 5} & {\bf 6} & {\bf 7} & {\bf 8} & n & {\bf 4} & {\bf 5} & {\bf 6} & {\bf 7} & {\bf 8}\\
\hline
\hline
2 & 1 & 0 & 0 & 0 & 0 &25 & 0 & 12 & 85 & 107 & 51 &48 & 21 & 45 & 613 & 1574 & 2045\\
3 & 0 & 1 & 0 & 0 & 0 &26 & 5 & 16 & 119 & 126 & 109 &49 & 0 & 40 & 614 & 1751 & 1554\\
4 & 2 & 0 & 1 & 0 & 0 &27 & 0 & 21 & 105 & 142 & 78 &50 & 10 & 54 & 771 & 1874 & 2505\\
5 & 0 & 1 & 1 & 0 & 0 &28 & 8 & 18 & 134 & 179 & 144 &51 & 0 & 66 & 704 & 1963 & 1946\\
6 & 2 & 2 & 2 & 0 & 1 &29 & 0 & 16 & 135 & 198 & 106 &52 & 13 & 58 & 771 & 2247 & 3008\\
7 & 0 & 3 & 1 & 1 & 0 &30 & 8 & 24 & 187 & 216 & 218 &53 & 0 & 48 & 788 & 2419 & 2322\\
8 & 4 & 1 & 5 & 1 & 1 &31 & 0 & 32 & 149 & 257 & 150 &54 & 12 & 66 & 989 & 2511 & 3713\\
9 & 0 & 2 & 5 & 1 & 1 &32 & 12 & 24 & 189 & 304 & 274 &55 & 0 & 92 & 849 & 2735 & 2829\\
10 & 3 & 3 & 9 & 3 & 2 &33 & 0 & 18 & 197 & 329 & 212 &56 & 18 & 68 & 938 & 3041 & 4354\\
11 & 0 & 5 & 7 & 4 & 1 &34 & 6 & 26 & 251 & 382 & 382 &57 & 0 & 49 & 1005 & 3187 & 3418\\
12 & 5 & 3 & 14 & 5 & 5 &35 & 0 & 37 & 218 & 431 & 279 &58 & 9 & 71 & 1175 & 3453 & 5233\\
13 & 0 & 4 & 14 & 7 & 2 &36 & 13 & 23 & 278 & 483 & 499 &59 & 0 & 98 & 1038 & 3659 & 4063\\
14 & 3 & 7 & 23 & 9 & 8 &37 & 0 & 24 & 275 & 547 & 366 &60 & 22 & 70 & 1215 & 3954 & 6234\\
15 & 0 & 10 & 17 & 12 & 5 &38 & 6 & 38 & 354 & 601 & 650 &61 & 0 & 63 & 1193 & 4315 & 4784\\
16 & 7 & 6 & 28 & 18 & 12 &39 & 0 & 45 & 313 & 643 & 493 &62 & 9 & 96 & 1440 & 4526 & 7301\\
17 & 0 & 6 & 27 & 22 & 8 &40 & 15 & 37 & 361 & 764 & 815 &63 & 0 & 104 & 1328 & 4674 & 5740\\
18 & 5 & 7 & 44 & 25 & 25 &41 & 0 & 30 & 359 & 838 & 623 &64 & 21 & 92 & 1378 & 5248 & 8514\\
19 & 0 & 12 & 35 & 36 & 13 &42 & 10 & 33 & 472 & 889 & 1083 &65 & 0 & 74 & 1440 & 5600 & 6631\\
20 & 7 & 9 & 54 & 46 & 30 &43 & 0 & 52 & 405 & 998 & 800 &66 & 14 & 80 & 1751 & 5741 & 10103\\
21 & 0 & 8 & 57 & 48 & 23 &44 & 11 & 44 & 480 & 1134 & 1305 &67 & 0 & 122 & 1531 & 6159 & 7794\\
22 & 4 & 15 & 77 & 62 & 51 &45 & 0 & 34 & 511 & 1197 & 1020 &68 & 16 & 98 & 1675 & 6730 & 11572\\
23 & 0 & 20 & 59 & 76 & 33 &46 & 7 & 56 & 609 & 1324 & 1653 &69 & 0 & 72 & 1792 & 7005 & 9097\\
24 & 11 & 11 & 87 & 88 & 76 &47 & 0 & 69 & 519 & 1435 & 1261 &70 & 14 & 120 & 2066 & 7465 & 13428\\
\hline
\hline
\end{array}
\end{equation*}
}
\end{table}

Take an $i$-hedrite $G$ with $i\in \{5,6,7\}$.
Then, if $F$ is a face of size $2$,
we reduce it to a vertex by using the following reduction operation:
\begin{center}
\begin{picture}(0,0)%
\includegraphics{ReductionOperation.pstex}%
\end{picture}%
\setlength{\unitlength}{1973sp}%
\begingroup\makeatletter\ifx\SetFigFont\undefined%
\gdef\SetFigFont#1#2#3#4#5{%
  \reset@font\fontsize{#1}{#2pt}%
  \fontfamily{#3}\fontseries{#4}\fontshape{#5}%
  \selectfont}%
\fi\endgroup%
\begin{picture}(4374,1128)(589,95)
\put(1576,1064){\makebox(0,0)[lb]{\smash{{\SetFigFont{6}{7.2}{\rmdefault}{\mddefault}{\updefault}{\color[rgb]{0,0,0}$F_1$}%
}}}}
\put(1576,239){\makebox(0,0)[lb]{\smash{{\SetFigFont{6}{7.2}{\rmdefault}{\mddefault}{\updefault}{\color[rgb]{0,0,0}$F_2$}%
}}}}
\put(4426,1064){\makebox(0,0)[lb]{\smash{{\SetFigFont{6}{7.2}{\rmdefault}{\mddefault}{\updefault}{\color[rgb]{0,0,0}$G_1$}%
}}}}
\put(4426,164){\makebox(0,0)[lb]{\smash{{\SetFigFont{6}{7.2}{\rmdefault}{\mddefault}{\updefault}{\color[rgb]{0,0,0}$G_2$}%
}}}}
\put(2401,614){\makebox(0,0)[lb]{\smash{{\SetFigFont{6}{7.2}{\rmdefault}{\mddefault}{\updefault}{\color[rgb]{0,0,0}$v'$}%
}}}}
\put(826,614){\makebox(0,0)[lb]{\smash{{\SetFigFont{6}{7.2}{\rmdefault}{\mddefault}{\updefault}{\color[rgb]{0,0,0}$v$}%
}}}}
\put(4651,614){\makebox(0,0)[lb]{\smash{{\SetFigFont{6}{7.2}{\rmdefault}{\mddefault}{\updefault}{\color[rgb]{0,0,0}$w$}%
}}}}
\put(1576,614){\makebox(0,0)[lb]{\smash{{\SetFigFont{6}{7.2}{\rmdefault}{\mddefault}{\updefault}{\color[rgb]{0,0,0}$F$}%
}}}}
\end{picture}%
\par
\end{center}
and get a graph denoted by $Red_F(G)$.
During this operation the vertices $v$ and $v'$ are merged
into one vertex $w$ and the faces $F_1$ and $F_2$ are
changed into $G_1$ and $G_2$ with one edge less.
Thus, it is possible that $G_1$ and/or $G_2$ are themselves of size $2$.
We apply the reduction operation whenever, by doing it, the reduced graph 
is still an $i$-hedrite.
Eventually, since every application of the technique diminish the vertex-set
one obtains a graph, denoted by $Red_{\infty}(G)$ for which we cannot apply the reduction operation anymore.

We call a graph {\em unreducible} if we cannot apply to it any reduction
operation. Let  $G'$ be an unreducible graph.
If $G'$ has no faces of size $2$, then it is an $8$-hedrite,
i.e. an octahedrite. If $G'$ has a face $F$ of size $2$, then denote
by $e_1$, $e_2$ the two edges of $F$.
Since $G'$ is unreducible, $F$ is adjacent on $e_1$ or $e_2$, say $e_1$,
to another face of size $2$.

\begin{figure}
\begin{center}
\resizebox{5.0cm}{!}{\includegraphics{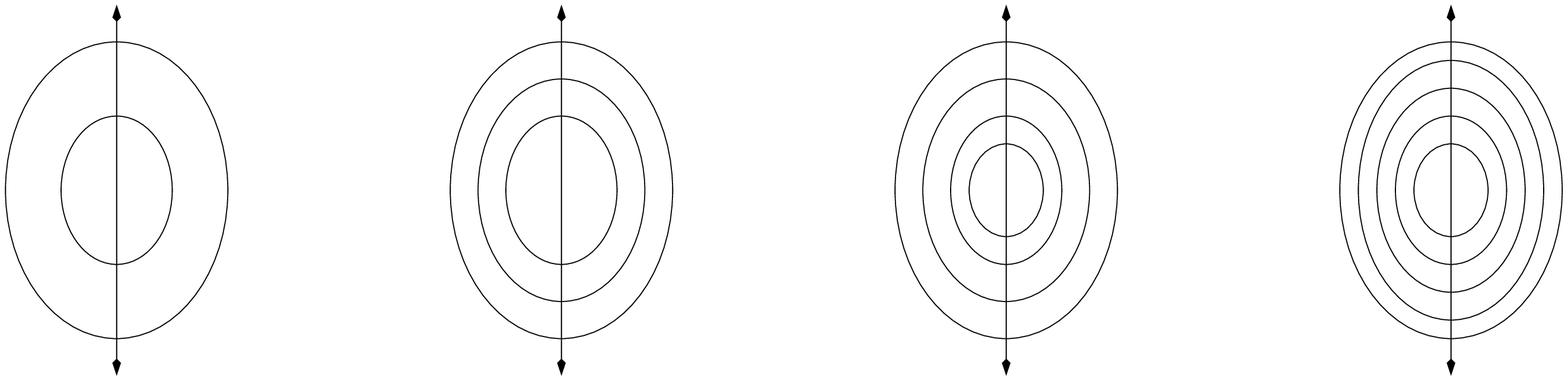}}\par
\end{center}
\caption{Infinite family of unreducible $4$-hedrites}
\label{InfiniteFamilyI4}
\end{figure}

If $e_2$ is incident to another face of size $2$, then $G'$ is actually
$2_{1}$, i.e. the unique graph with two vertices, and four faces of size
$2$, i.e. the 1st one on Figure  \ref{InfiniteFamilyI4}.
It is easy to see that $e_2$ cannot be incident to another face of size $3$,
but it can be incident to another face of size $4$ and in that case
$G'$ is not $3$-connected and thus (see \cite{oct2})
it belongs to the infinite family depicted in
Figure \ref{InfiniteFamilyI4}.

Call {\em expansion operation} the reverse of the reduction operation.
The generation method of $i$-hedrite is to consider all unreducible
$i$-hedrites and all possible ways of expanding them.
For an unreducible graph $G$ denote by ${\mathcal Exp}(G)$ the set of all
possible $i$-hedrites that can be obtained by repeated application of the
expansion operation.
For the graphs of the infinite family of Figure \ref{InfiniteFamilyI4}
no expansion operation is possible and thus no $i$-hedrite is obtained
from them.
A priori, the set ${\mathcal Exp}(G)$ can be infinite 
but, as far as we know, for any $8$-hedrite $G$ the set
${\mathcal Exp}(G)$ is finite although we have no proof
of it.
It turns out that ${\mathcal Exp}(2_1)$ is infinite but it has a simple
description.

\begin{figure}
\begin{center}
\begin{minipage}[b]{2.3cm}
\centering
\resizebox{16mm}{!}{\rotatebox{0}{\includegraphics{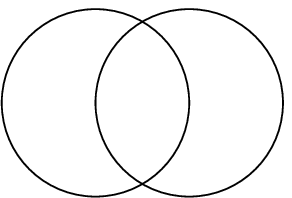}}}\par
$D_{4h}$, $n=2$
\end{minipage}
\begin{minipage}[b]{2.3cm}
\centering
\resizebox{16mm}{!}{\rotatebox{0}{\includegraphics{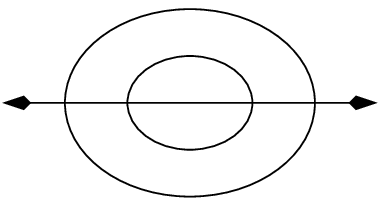}}}\par
$D_{2h}$, $n=4$
\end{minipage}
\begin{minipage}[b]{2.3cm}
\centering
\resizebox{16mm}{!}{\rotatebox{0}{\includegraphics{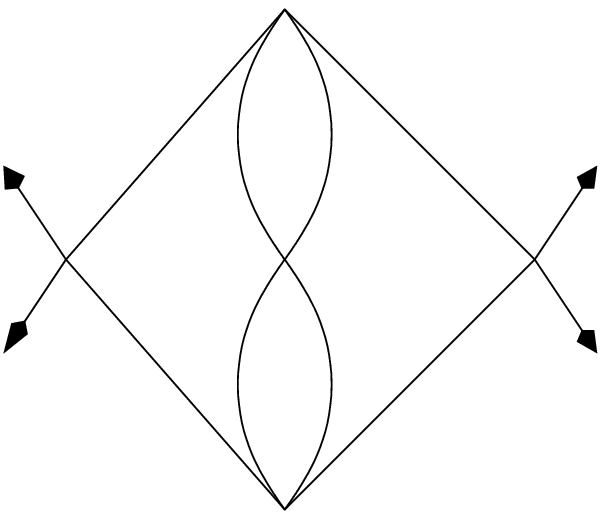}}}\par
$D_{2d}$, $n=6$
\end{minipage}
\begin{minipage}[b]{2.3cm}
\centering
\resizebox{16mm}{!}{\rotatebox{0}{\includegraphics{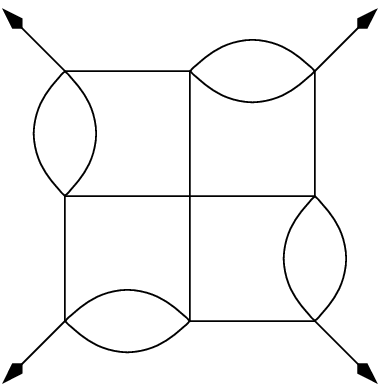}}}\par
$D_{4}$, $n=10$
\end{minipage}
\begin{minipage}[b]{2.3cm}
\centering
\resizebox{16mm}{!}{\rotatebox{0}{\includegraphics{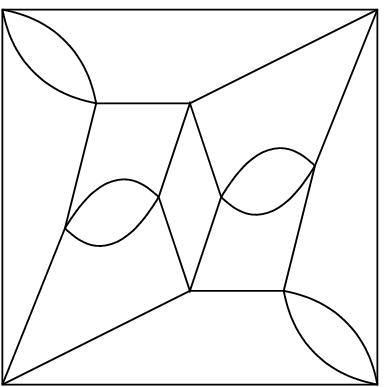}}}\par
$D_{2}$, $n=12$
\end{minipage}

\end{center}
\caption{Minimal representatives for each possible symmetry group of a $4$-hedrite}
\label{MinimalRepresentative4hedrite}
\end{figure}

\begin{figure}
\begin{center}
\begin{minipage}[b]{2.3cm}
\centering
\resizebox{16mm}{!}{\rotatebox{0}{\includegraphics{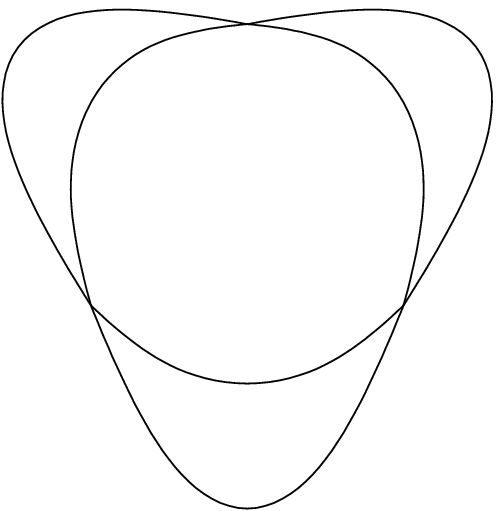}}}\par
$D_{3h}$, $n=3$
\end{minipage}
\begin{minipage}[b]{2.3cm}
\centering
\resizebox{16mm}{!}{\rotatebox{0}{\includegraphics{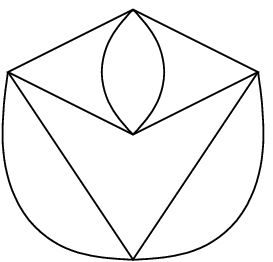}}}\par
$C_{2v}$, $n=5$
\end{minipage}
\begin{minipage}[b]{2.3cm}
\centering
\resizebox{16mm}{!}{\rotatebox{0}{\includegraphics{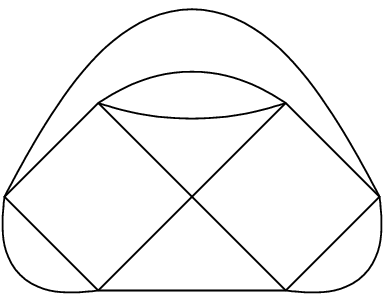}}}\par
$C_{s}$, $n=7$
\end{minipage}
\begin{minipage}[b]{2.3cm}
\centering
\resizebox{16mm}{!}{\rotatebox{0}{\includegraphics{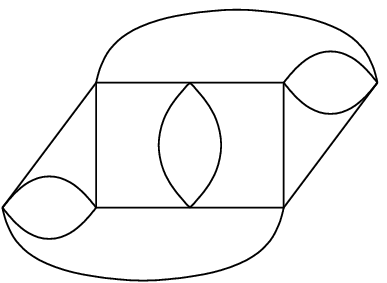}}}\par
$C_{2}$, $n=8$
\end{minipage}
\begin{minipage}[b]{2.3cm}
\centering
\resizebox{16mm}{!}{\rotatebox{0}{\includegraphics{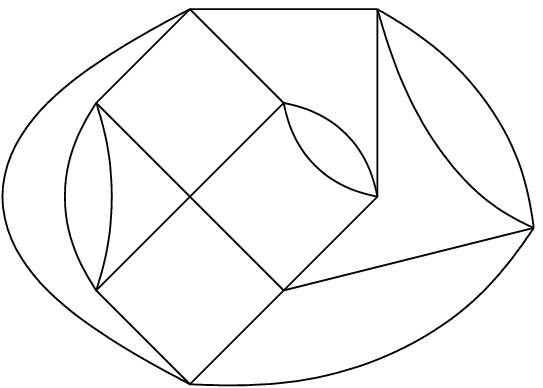}}}\par
$C_{1}$, $n=10$
\end{minipage}
\begin{minipage}[b]{2.3cm}
\centering
\resizebox{16mm}{!}{\rotatebox{0}{\includegraphics{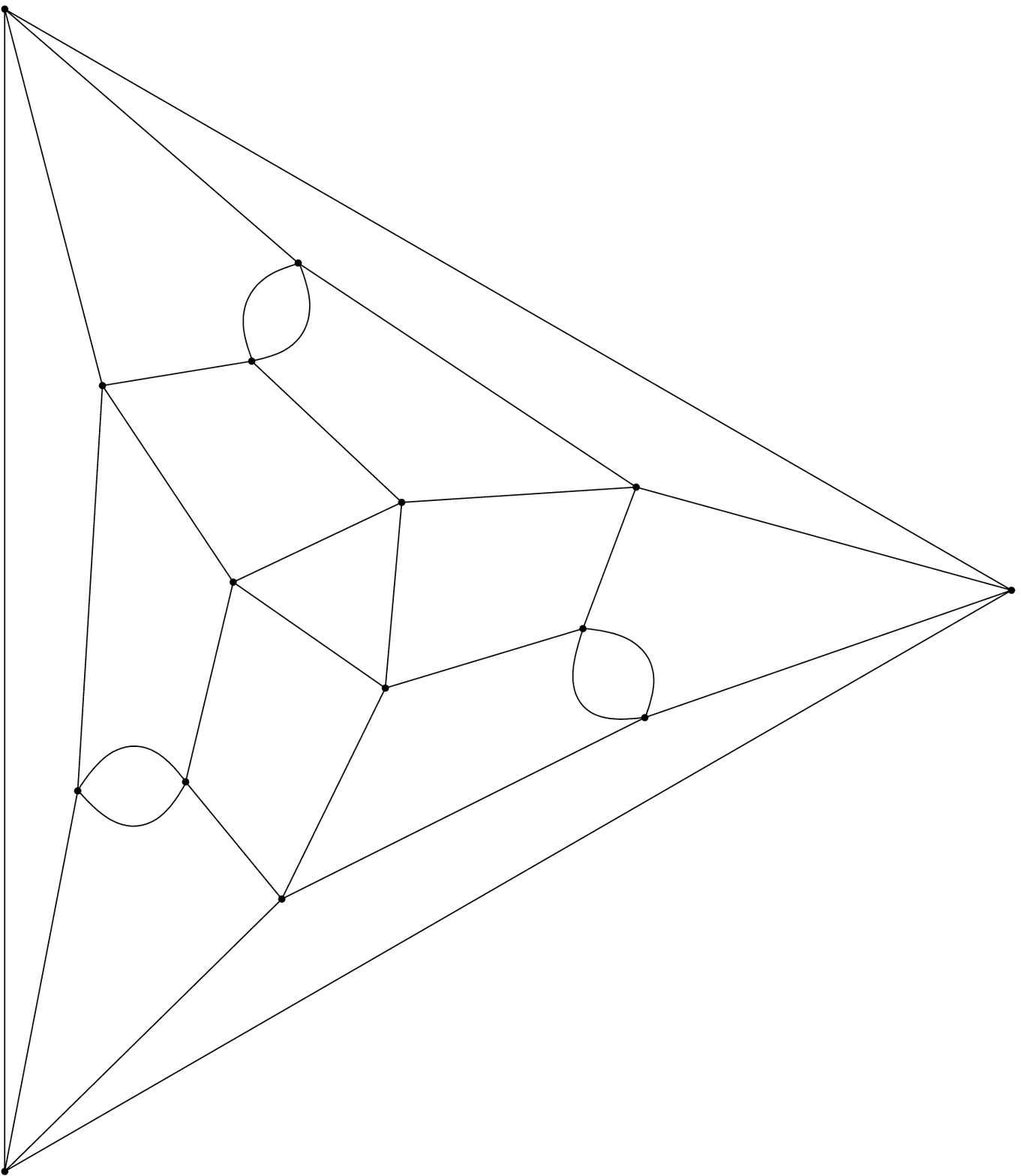}}}\par
$D_{3}$, $n=15$
\end{minipage}

\end{center}
\caption{Minimal representatives for each possible symmetry group of a $5$-hedrite}
\label{MinimalRepresentative5hedrite}
\end{figure}

\begin{figure}
\begin{center}
\begin{minipage}[b]{2.3cm}
\centering
\resizebox{16mm}{!}{\rotatebox{0}{\includegraphics{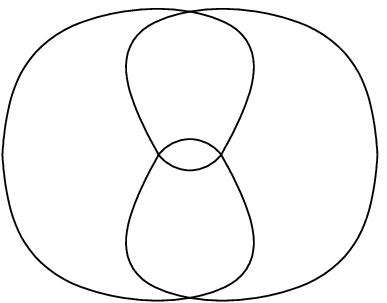}}}\par
$D_{2d}$, $n=4$
\end{minipage}
\begin{minipage}[b]{2.3cm}
\centering
\resizebox{16mm}{!}{\rotatebox{0}{\includegraphics{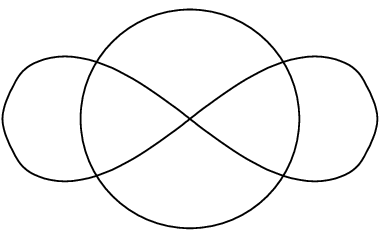}}}\par
$C_{2v}$, $n=5$
\end{minipage}
\begin{minipage}[b]{2.3cm}
\centering
\resizebox{16mm}{!}{\rotatebox{0}{\includegraphics{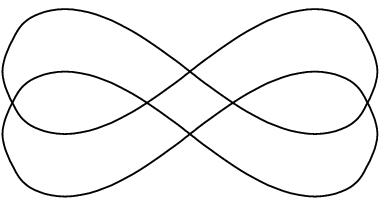}}}\par
$D_{2h}$, $n=6$
\end{minipage}
\begin{minipage}[b]{2.3cm}
\centering
\resizebox{16mm}{!}{\rotatebox{0}{\includegraphics{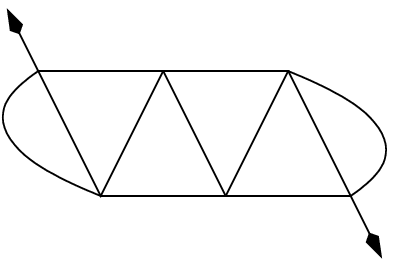}}}\par
$C_2$, $n=6$
\end{minipage}
\begin{minipage}[b]{2.3cm}
\centering
\resizebox{16mm}{!}{\rotatebox{0}{\includegraphics{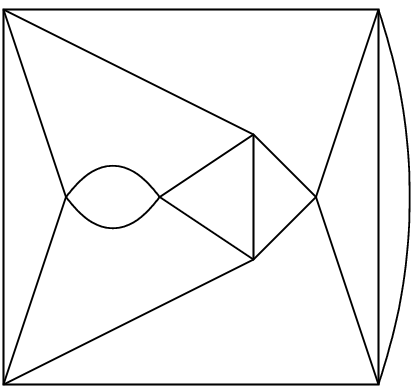}}}\par
$C_s$, $n=9$
\end{minipage}
\begin{minipage}[b]{2.3cm}
\centering
\resizebox{16mm}{!}{\rotatebox{0}{\includegraphics{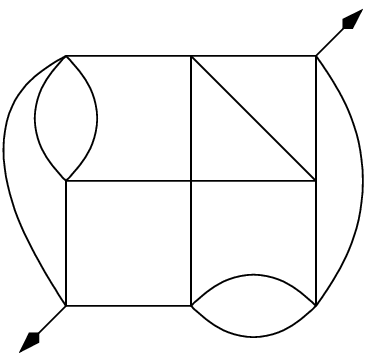}}}\par
$C_1$, $n=9$
\end{minipage}
\begin{minipage}[b]{2.3cm}
\centering
\resizebox{16mm}{!}{\rotatebox{0}{\includegraphics{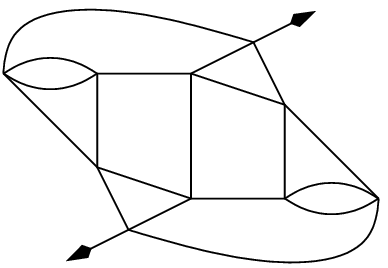}}}\par
$C_{2h}$, $n=10$
\end{minipage}
\begin{minipage}[b]{2.3cm}
\centering
\resizebox{16mm}{!}{\rotatebox{0}{\includegraphics{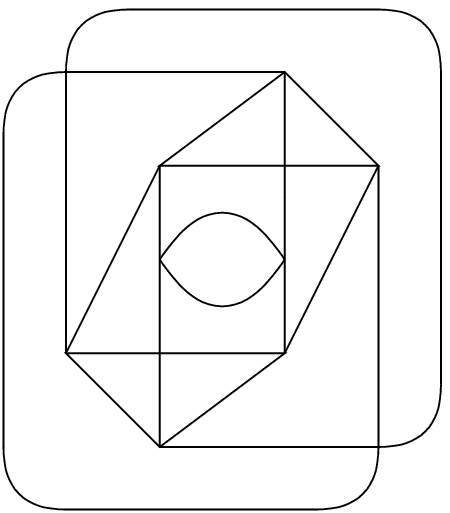}}}\par
$D_2$, $n=12$
\end{minipage}
\begin{minipage}[b]{2.3cm}
\centering
\resizebox{16mm}{!}{\rotatebox{0}{\includegraphics{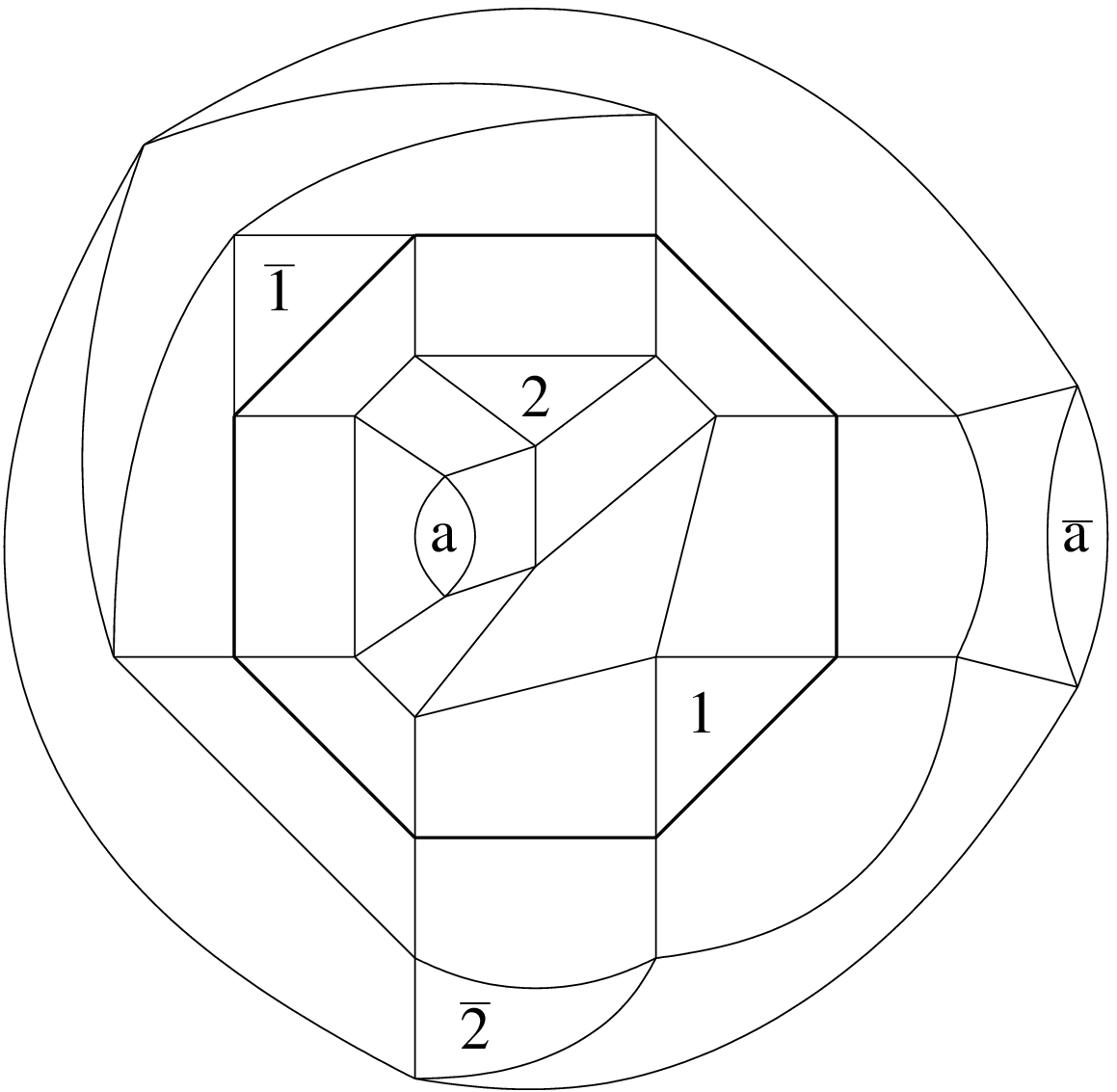}}}\par
$C_i$, $n=30$
\end{minipage}

\end{center}
\caption{Minimal representatives for each possible symmetry group of a $6$-hedrite}
\label{MinimalRepresentative6hedrite}
\end{figure}

\begin{figure}
\begin{center}
\begin{minipage}[b]{2.3cm}
\centering
\resizebox{16mm}{!}{\rotatebox{0}{\includegraphics{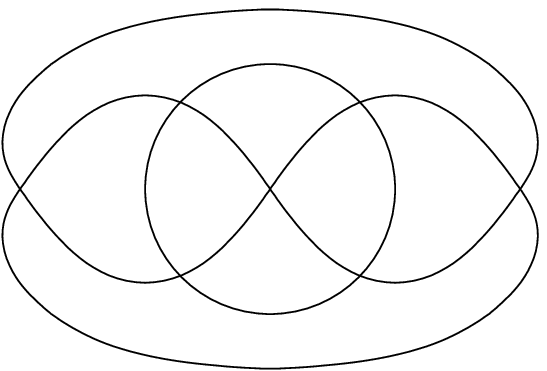}}}\par
$C_{2v}$, $n=7$
\end{minipage}
\begin{minipage}[b]{2.3cm}
\centering
\resizebox{16mm}{!}{\rotatebox{0}{\includegraphics{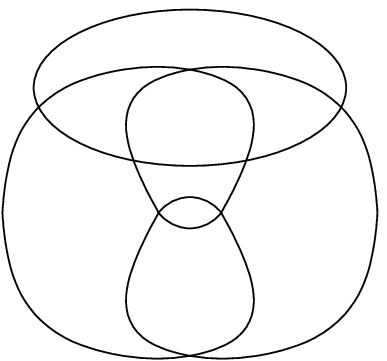}}}\par
$C_s$, $n=8$
\end{minipage}
\begin{minipage}[b]{2.3cm}
\centering
\resizebox{16mm}{!}{\rotatebox{0}{\includegraphics{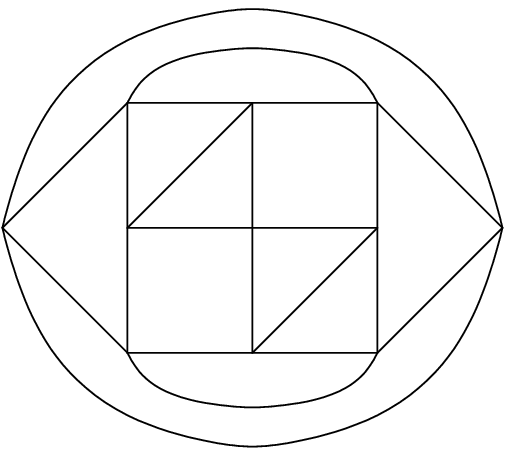}}}\par
$C_2$, $n=11$
\end{minipage}
\begin{minipage}[b]{2.3cm}
\centering
\resizebox{16mm}{!}{\rotatebox{0}{\includegraphics{7-hedrite11_1sec.eps}}}\par
$C_1$, $n=11$
\end{minipage}

\end{center}
\caption{Minimal representatives for each possible symmetry group of a $7$-hedrite}
\label{MinimalRepresentative7hedrite}
\end{figure}

\begin{theorem}
(i) The only symmetry groups of $4$-hedrites are $D_{4h}$, $D_4$, $D_{2h}$, $D_{2d}$ and $D_2$.

(ii) The only symmetry groups of $5$-hedrites are: 
$D_{3h}$, $D_3$, $C_{2v}$, $C_s$, $C_2$ and $C_1$.

(iii) The only symmetry groups of $6$-hedrites are: $D_{2d}$, $D_{2h}$, $D_{2}$, $C_{2h}$, $C_{2v}$, $C_i$, $C_{2}$, $C_{s}$, $C_{1}$.

(iv) The only symmetry groups of $7$-hedrites are: $C_{2v}$, $C_{2}$, $C_{s}$ and $C_{1}$.
\end{theorem}
The theorem is proven in the same way as for octahedrites. Minimal representative for each symmetry groups are given in Figures \ref{MinimalRepresentative4hedrite}, \ref{MinimalRepresentative5hedrite}, \ref{MinimalRepresentative6hedrite}, \ref{MinimalRepresentative7hedrite}.

Further generalization of octahedrites are $4$-regular plane graphs with
$4$-, $3$-, $2$- and $1$-gonal faces only. Then, besides $i$-hedrites, we get
graphs with $(p_1,p_2,p_3)$=$(2,1,0)$, $(2,0,2)$, $(1,2,1)$, $(1,1,3)$,
$(1,0,5)$.
The enumeration method is then to use $i$-hedrites and to add a $1$-gon
when we have a pair of $2$-gon and $3$-gon that are adjacent in all possible
ways. This is simlar to the strategy of squeezing of $2$-gons used for
the enumeration of $i$-hedrites.

\section{Central circuits and alternating knots}\label{CentralCircuit}
The edges of an octahedrite, as of any Eulerian plane graph, are partitioned
by its {\em central circuits}, i.e. those which are obtained by starting
with an edge and continuing at each vertex by the edge opposite
the entering one. 
The central circuits of an octahedrite can define circle in the plane or
have self-intersections.

If $C_1$, $C_2$ are two (possibly, self-intersecting) central circuits
of an octahedrite $G$, then they
are called {\em parallel} if they are separated by a sequence of faces
of size $4$
(such pair is called {\em railroad} in \cite{oct2,octa}).
It is possible to reduce those two central circuits into just one and
thus get an octahedrite with less vertices.
We call an octahedrite {\em irreducible} if it has no parallel central
circuits.
Of course, the reverse operation is possible, i.e. split a central
circuit into two or more parallel central circuits. In this way
every octahedrite is obtained from an irreducible octahedrite.

It is proved in \cite{octa} that an irreducible octahedrite has at most
$6$ central circuits and in \cite{oct2} that an irreducible $i$-hedrite
has at most $i-2$ central circuits.
All irreducible octahedrites with non self-intersecting central
circuits have been classified in \cite{octa} (see, for another presentation,
\cite{oct2}).

\begin{figure}
\begin{center}
\begin{minipage}[b]{3.0cm}
\centering
\resizebox{2.8cm}{!}{\includegraphics{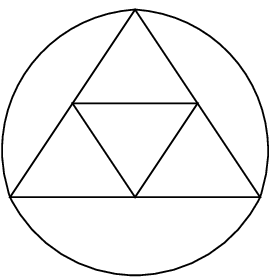}}\par
$O_{h}$, $n=6$
\end{minipage}
\begin{minipage}[b]{3.0cm}
\centering
\resizebox{2.8cm}{!}{\includegraphics{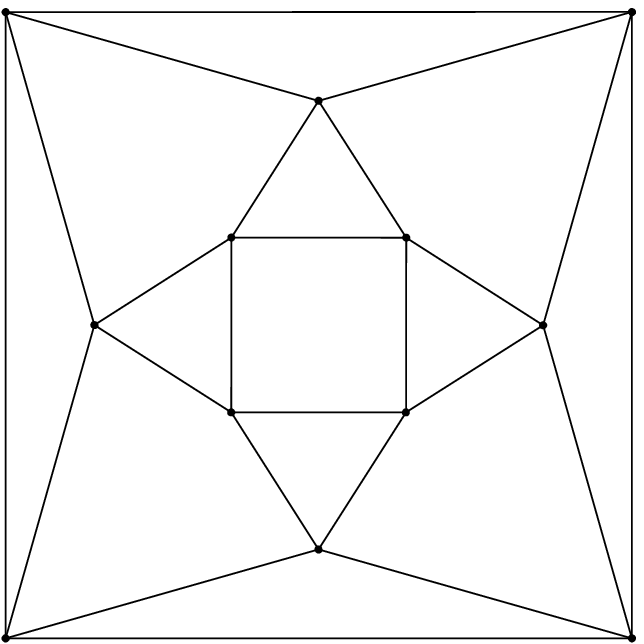}}\par
$O_{h}$, $n=12$
\end{minipage}
\begin{minipage}[b]{3.0cm}
\centering
\resizebox{2.8cm}{!}{\includegraphics{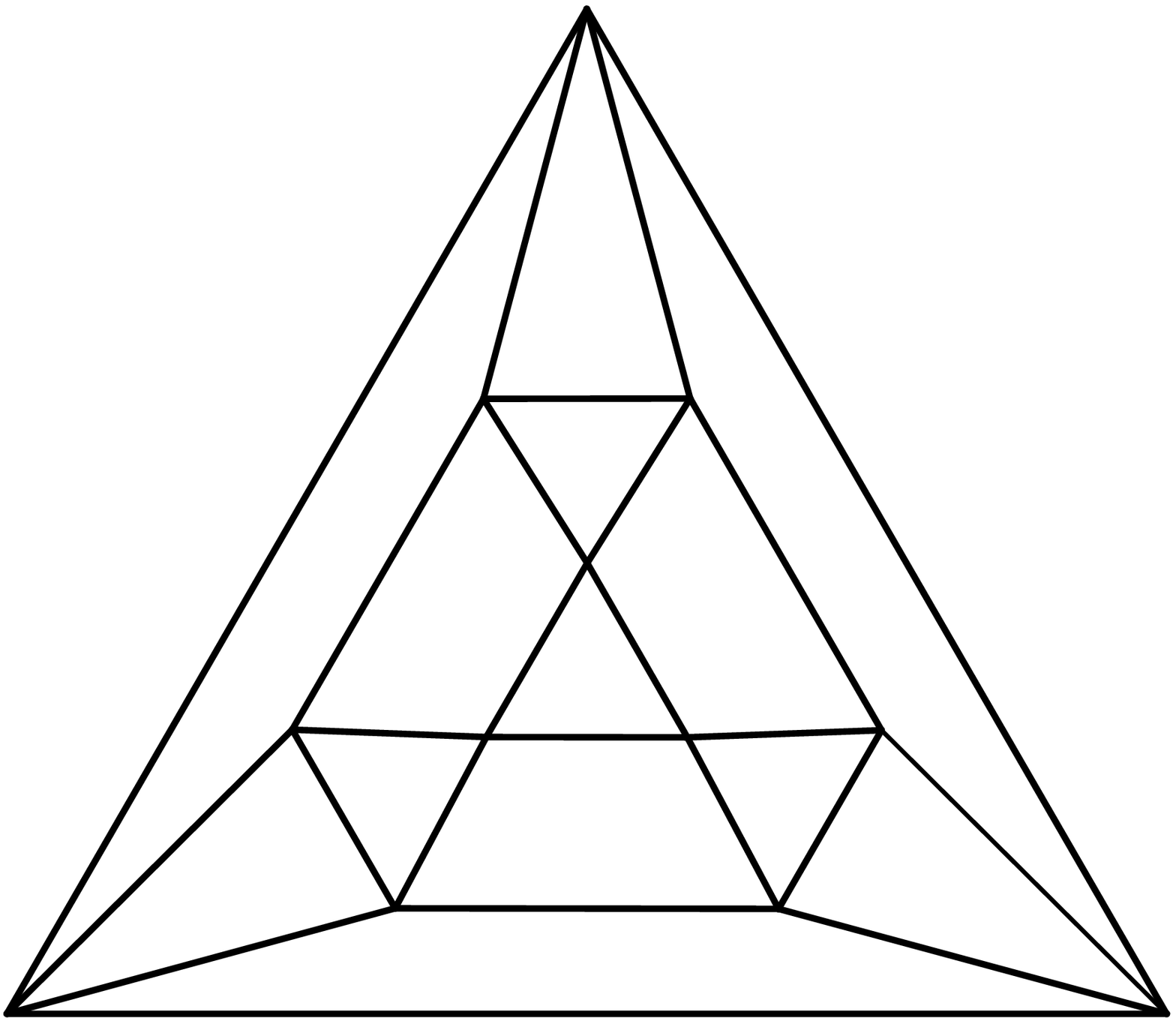}}\par
$D_{3h}$, $n=12$
\end{minipage}
\begin{minipage}[b]{3.0cm}
\centering
\resizebox{2.8cm}{!}{\includegraphics{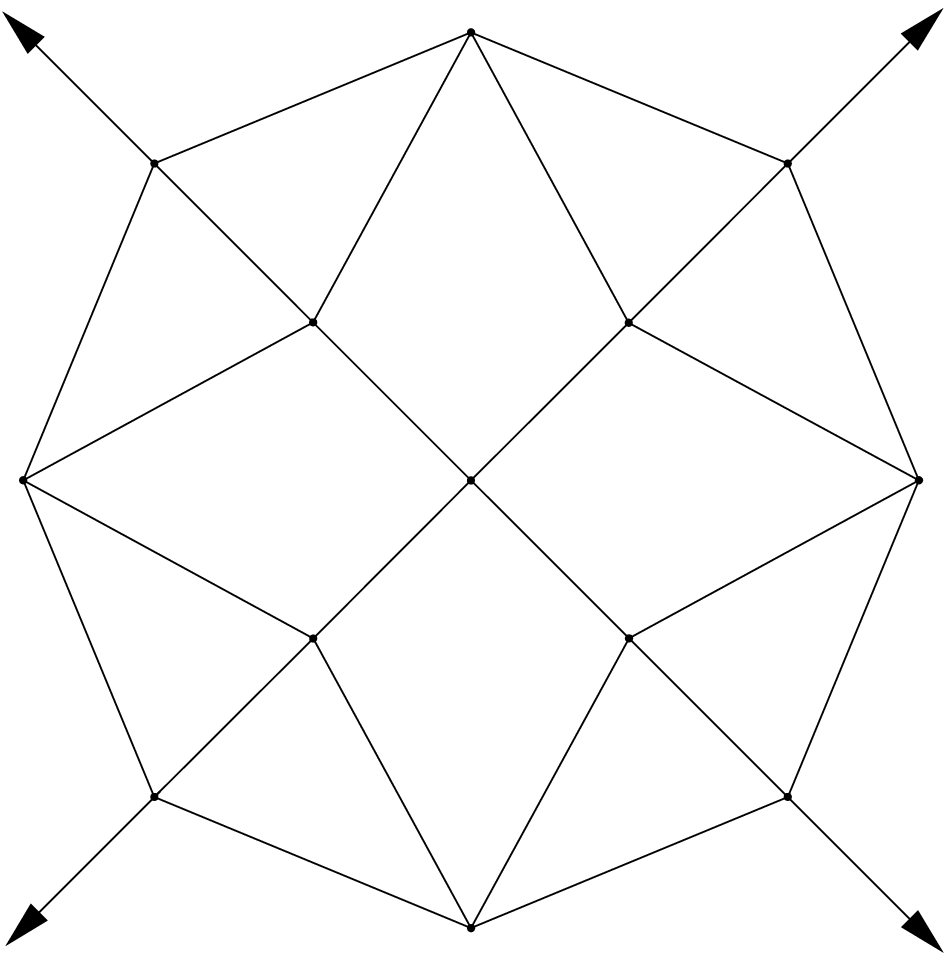}}\par
$D_{4h}$, $n=14$
\end{minipage}
\begin{minipage}[b]{3.0cm}
\centering
\resizebox{2.8cm}{!}{\includegraphics{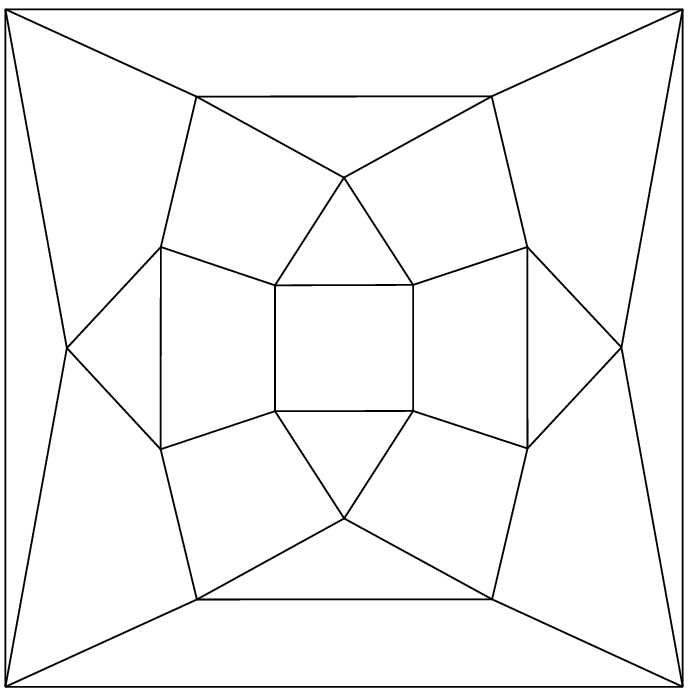}}\par
$D_{2d}$, $n=20$
\end{minipage}
\begin{minipage}[b]{3.0cm}
\centering
\resizebox{2.8cm}{!}{\includegraphics{PL2_22_D2hsec.eps}}\par
$D_{2h}$, $n=22$
\end{minipage}
\begin{minipage}[b]{3.0cm}
\centering
\resizebox{2.8cm}{!}{\includegraphics{PL2_30_Osec.eps}}\par
$O$, $n=30$
\end{minipage}
\begin{minipage}[b]{3.0cm}
\centering
\resizebox{2.8cm}{!}{\includegraphics{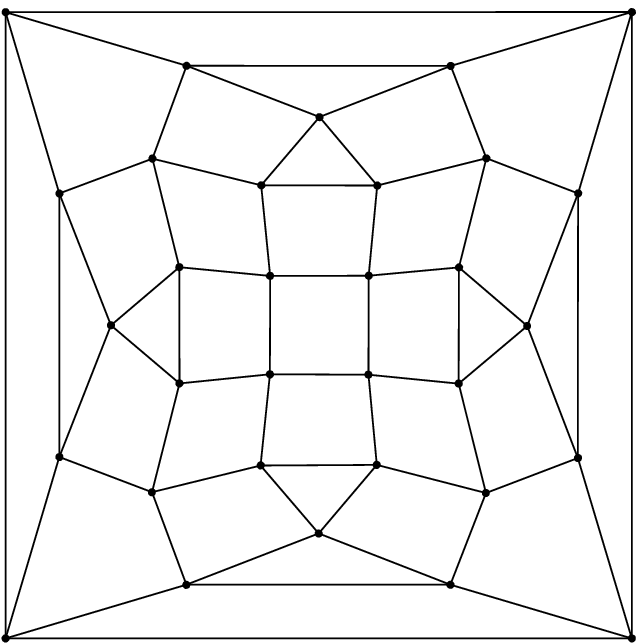}}\par
$D_{4h}$, $n=32$
\end{minipage}
\end{center}

\caption{The irreducible octahedrites with simple central circuits}
\label{IrreducibleSimpleOcta}
\end{figure}

\begin{theorem}
There are exactly eight irreducible octahedrites with simple central circuits (see Figure \ref{IrreducibleSimpleOcta}).
\end{theorem}

\begin{figure}
\begin{center}
\resizebox{7.2cm}{!}{\includegraphics{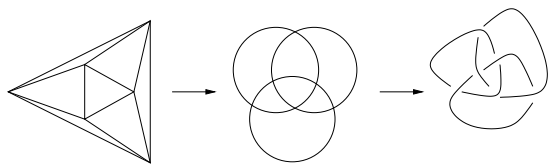}}\par
\end{center}
\caption{The link corresponding to the Octahedron}
\label{TheLinkOctahedron}
\end{figure}

A {\em link} is a set of circles embedded in $3$-space that do not intersect;
a link can be represented with its overlapping and underlapping on the
plane. A link with only one component is called a {\em knot} and Knot Theory
is concerned with characterizing different plane presentations
of links (see \cite{lickorish} for a pleasant introduction).
A link is called {\em alternating} if it admits a plane
representation in which overlappings and underlappings alternate.
For a $4$-regular plane graph we can define a corresponding alternating
link, where the central circuits correspond to the components of the link
(see an example on Figure \ref{TheLinkOctahedron}).
It is interesting that there is no known topological characterization of 
alternating links.

\begin{figure}
\begin{center}
\begin{minipage}[b]{3.4cm}
\centering
\resizebox{3.2cm}{!}{\includegraphics{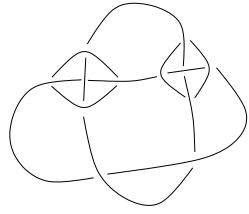}}\par
\end{minipage}

\end{center}
\caption{A Borromean link}
\label{TheCountExample}
\end{figure}

Since an octahedrite with $n$ vertices is $3$-connected, there is no
disjointing vertex and thus \cite[Chapter 5]{lickorish}) the corresponding
alternating link cannot be represented with less than 
$n$ crossings. But it can happen that two octahedrites that are not 
equivalent as graphs give rise to equivalent alternating links.
A link with $m$ components is called {\em Borromean} if after removal
of any $m-2$ components the remaining two components can be
separated one from the other.
It is conjectured in \cite{octa} that an alternating link obtained
from a $4$-regular
$3$-connected plane graph is Borromean if and only if for any two
central circuits the distance between any two of its consecutive points
of its intersection is even.
This condition is, of course, 
sufficient but there are reasons to think
that it is not necessary since there exist $4$-regular plane graphs 
(but not $3$-connected) which are 
Borromean
without satisfying the specified condition, see Figure \ref{TheCountExample}.

\section{Self-dual graphs}
A graph $G$ is called {\em self-dual} if it is isomorphic to its dual $G^{*}$.
The {\em medial graph} $Med(G)$ of a plane graph $G$ is the plane graph
obtained by putting a vertex on any edge with two edges adjacent if they
share a common vertex and are contained in a common face.
One has $Med(G)=Med(G^*)$.
The graph $G'=Med(G)$ is always $4$-regular and its dual $(Med(G))^{*}$ is
bipartite, that is the face-set ${\mathcal F}$ of $Med(G)$ is split into
two sets ${\mathcal F}_1(G')$ and ${\mathcal F}_2(G')$, which correspond to
the vertices and faces of the graph $G$.
The bipartition ${\mathcal F}_1(G')$, ${\mathcal F}_2(G')$ can be computed
easily from a given $4$-regular plane graph, i.e. one can compute
easily from a graph $G'$ the two dual graphs $G_1$ and $G_2$
such that $G'=Med(G_1)=Med(G_2)$.

Call $G$ a 
{\em $i$-self-hedrite} if it is a self-dual plane graph
with vertices of degree $2$, $3$ or $4$ with $v_2+v_3=i$
and, consequently, faces of size $2$, $3$ or $4$. 
If $G$ is a $i$-self-hedrite then $Med(G)$ is a $2i$-hedrite.

The Euler formula $V-E+F=2$ for a self-dual plane graph is, clearly:
\begin{equation}
\sum_{j=2}^{\infty}p_j(4-j)=4;
\end{equation}
we again permit $2$-gons but not $1$-gons.
Define an {\em $i$-self-hedrite} to be a such graph with
faces of size $2, 3, 4$ only and $p_2+p_3=i$.
So, $2p_2+p_3=4$ and $p_4$ is not bounded; also $n=p_4+\frac{p_3}{2}+2=p_4-p_2+4$.
Clearly, an $i$-self-hedrite can have $i=2, 3, 4$ only with $(p_2,p_3)=(2,0), (1,2), (0,4)$, 
respectively.
The $i$-self-hedrites with smallest number $n$ of vertices have no $4$-gons; they are $Bundle_2$ ($2$ vertices connected by $2$ edges), 
triangle with one doubled edge and Tetrahedron, respectively.
Easy to check that $n$-vertex an $i$-self-hedrite exists if $n\ge i$.

Thus our enumeration method for $i$-self-hedrites is to consider all 
$2i$-hedrites $G'$,
determine for them the graphs $G_1$, $G_2$ such that $G'=Med(G_1)=Med(G_2)$
and keep the ones that have $G_1$ isomorphic to $G_2$.
We denote by $Med^{-1}(G')=G_1\simeq G_2$ the obtained plane graph if it exists.
Using the enumeration of $2i$-hedrites, we can derive the $i$-self-hedrite,
see Table \ref{TableSelfDual}. Another method would be possible
with the results of \cite{Archdeacon} but it would require more 
hard programming work and the speed gain is uncertain.

\begin{table}
\caption{Number of $i$-self-hedrites with $4\leq n\leq 40$ and $2\leq i\leq 4$}
\label{TableSelfDual}
\begin{equation*}
{\scriptsize
\begin{array}{||c|c|c|c||c|c|c|c||c|c|c|c||c|c|c|c||}
\hline
\hline
n & {\bf 2} & {\bf 3} & {\bf 4} & n & {\bf 2} & {\bf 3} & {\bf 4} & n & {\bf 2} 
& {\bf 3} & {\bf 4} & n & {\bf 2} & {\bf 3} & {\bf 4}\\
\hline
\hline
2 & 1 & 0 & 0 &12 & 4 & 29 & 24 &22 & 10 & 90 & 191 &32 & 9 & 239 & 584\\
3 & 1 & 1 & 0 &13 & 6 & 30 & 33 &23 & 7 & 119 & 198 &33 & 9 & 256 & 631\\
4 & 2 & 1 & 1 &14 & 5 & 42 & 40 &24 & 7 & 131 & 234 &34 & 14 & 232 & 748\\
5 & 2 & 4 & 1 &15 & 5 & 47 & 48 &25 & 10 & 124 & 276 &35 & 10 & 290 & 760\\
6 & 3 & 6 & 2 &16 & 8 & 48 & 69 &26 & 10 & 162 & 304 &36 & 14 & 308 & 857\\
7 & 3 & 7 & 4 &17 & 5 & 64 & 73 &27 & 8 & 170 & 332 &37 & 16 & 286 & 956\\
8 & 3 & 11 & 6 &18 & 6 & 72 & 92 &28 & 12 & 158 & 407 &38 & 11 & 342 & 1002\\
9 & 3 & 16 & 8 &19 & 8 & 70 & 114 &29 & 10 & 190 & 421 &39 & 11 & 359 & 1070\\
10 & 5 & 16 & 15 &20 & 6 & 89 & 130 &30 & 9 & 210 & 476 &40 & 16 & 332 & 1239\\
11 & 4 & 26 & 16 &21 & 8 & 104 & 148 &31 & 14 & 202 & 550 & & & & \\
\hline
\hline
\end{array}
}
\end{equation*}
\end{table}

\begin{theorem}
(i) The possible symmetry groups of $2$-self-hedrites graphs are $C_2$, $C_{2v}$, $C_{2h}$, $D_{2}$ and $D_{2h}$. Minimal representatives are given in Figure \ref{MinimalRepresentative2SelfHedrites}.

(ii) The possible symmetry groups of $3$-self-hedrites graphs are $C_1$, $C_2$, $C_s$ and $C_{2v}$. Minimal representatives are given in Figure \ref{MinimalRepresentative3SelfHedrites}.

(iii) The possible symmetry groups of $4$-self-hedrites graphs are $C_1$, $C_2$, $C_{2h}$, $C_{2\nu}$, $C_{3}$, $C_{3\nu}$, $C_{4}$, $C_{4\nu}$, $C_{i}$, $C_{s}$, $D_{2}$, $D_{2d}$, $D_{2h}$, $S_{4}$, $T$, $T_d$.
Minimal representatives are given in Figure \ref{MinimalRepresentativeSelfDual34}.

\end{theorem}
\proof If $G$ is a $4$-self-hedrite then $G'=Med(G)$ is an octahedrite.
If $\Gamma$, $\Gamma'$ are the symmetry groups of $G$, $G'$, then the
self-duality of $G$ becomes a symmetry in $G'$ that exchanges ${\mathcal F}_1(G')$ and ${\mathcal F}_2(G')$.
Thus $\Gamma$ is identified with the subgroup of $\Gamma'$ formed by
the transformations preserving the bipartition. Obviously, the order
of $\Gamma$ is half the one of $\Gamma'$.
The possible groups of $G'$ are known (see Theorem \ref{GroupsOfOctahedrites}).
So, we set out to enumerate the index $2$ subgroups of each of the $18$
groups and 
found, besides the groups in the statement, the groups
$D_3$, $D_4$, $C_{4h}$, $C_{3h}$, $S_6$, $S_8$ and $T_h$.

The graph $G$ has $4$ vertices of degree $3$ and $4$ faces of size 
$3$; both should be partitioned in the same number of orbits and this excludes
$D_4$, $D_3$, $C_{3h}$, $S_6$ and $T_h$.
Suppose $G$ has symmetry $C_{4h}$. Due to the plane of 
symmetry, the $4$-fold axis pass 
through, either two vertices of degree $4$, or through 
two faces of size $4$. But self-duality requires that it passes through
a vertex and a face. The same argument excludes $S_8$.

For $2$-self hedrites, using the known groups for $4$-hedrites gives 
candidates $C_2$, $C_{2h}$, $C_{2v}$, $C_4$, $C_{4h}$, $C_{4v}$, $D_2$, 
$D_{2d}$, $D_{2h}$, $D_4$, $S_4$. Same kind of orbit reasons exclude
$C_4$, $C_{4h}$, $C_{4v}$, $D_{2d}$, $D_4$, $S_4$.
A $3$-self-hedrite has only one vertex of degree $2$ that has to be
preserved by any symmetry. So, the symmetry is a subgroup of $C_{2v}$
and all possible subgroups do occur. \qed

\begin{figure}
\begin{center}
\begin{minipage}[b]{3.0cm}
\centering
\resizebox{2.9cm}{!}{\includegraphics{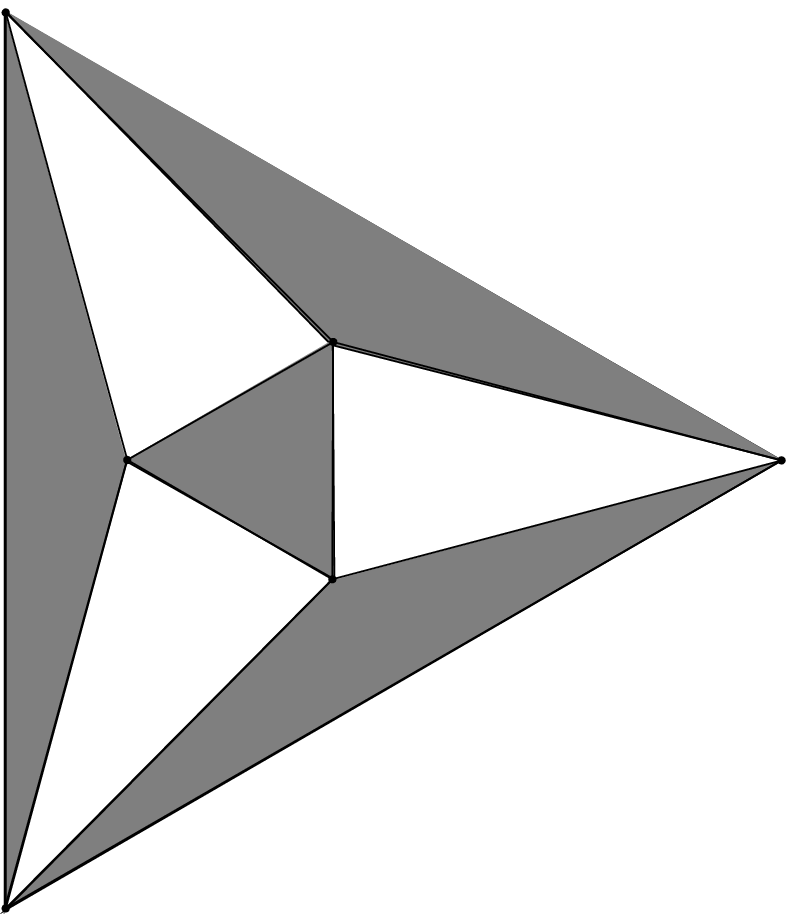}}\par
$(k,l)=(1,0)$
\end{minipage}
\begin{minipage}[b]{3.0cm}
\centering
\resizebox{2.9cm}{!}{\includegraphics{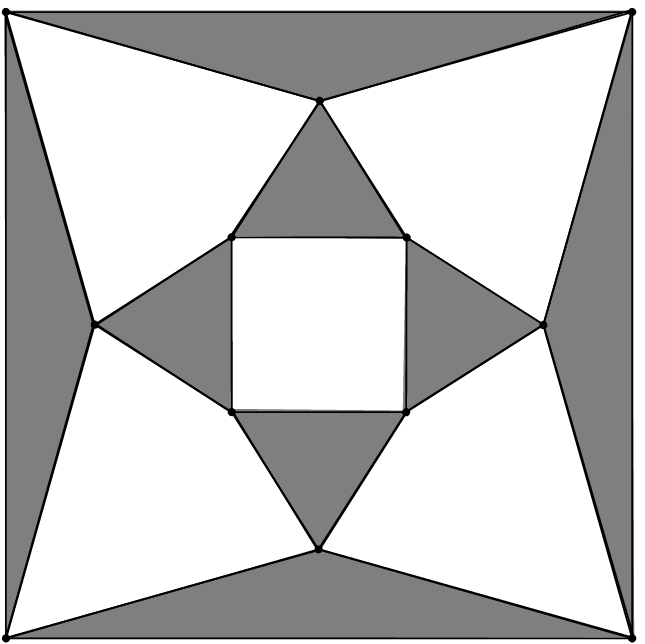}}\par
$(k,l)=(1,1)$
\end{minipage}
\begin{minipage}[b]{3.0cm}
\centering
\resizebox{2.9cm}{!}{\includegraphics{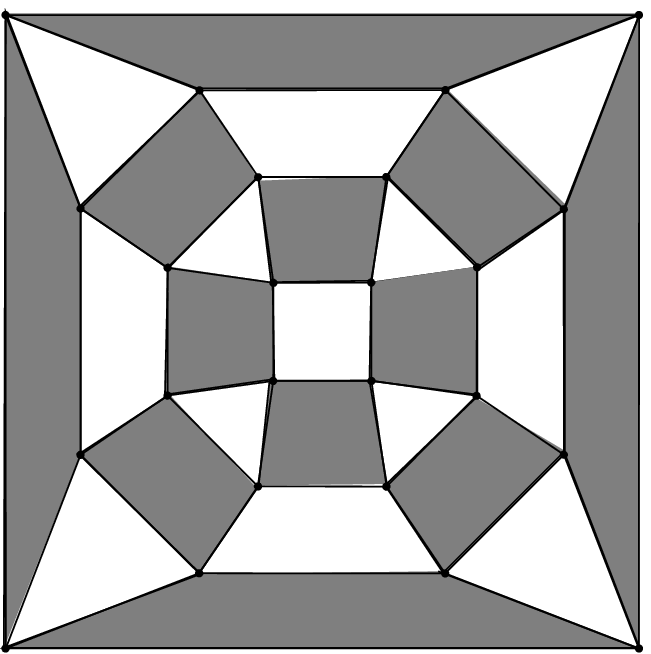}}\par
$(k,l)=(2,0)$
\end{minipage}
\begin{minipage}[b]{3.0cm}
\centering
\resizebox{2.9cm}{!}{\includegraphics{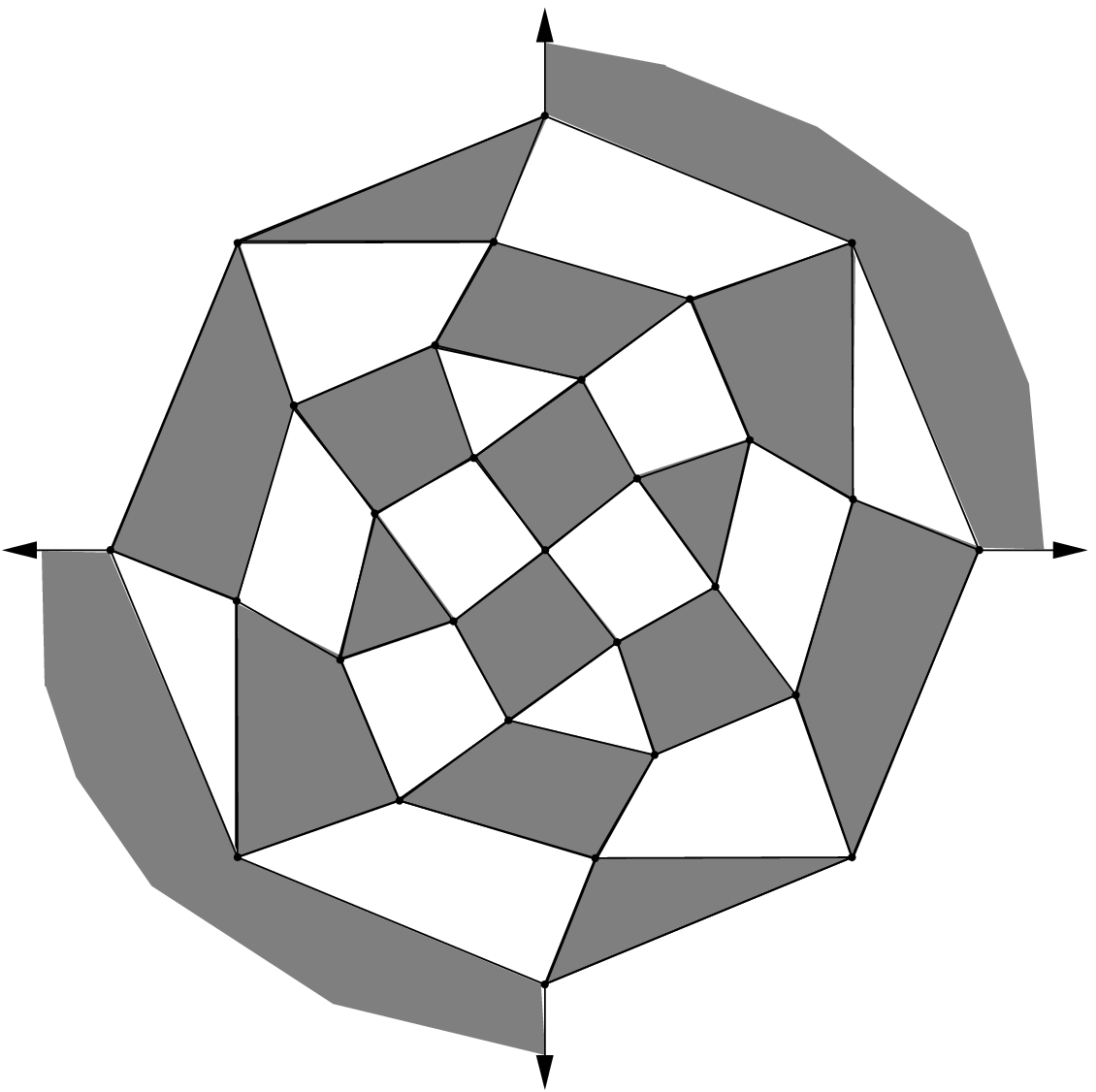}}\par
$(k,l)=(2,1)$
\end{minipage}

\end{center}
\caption{First examples of octahedrites of symmetry $O$ or $O_h$ expressed as $GC_{k,l}(Octahedron)$}
\label{FirstExample}
\end{figure}

It is known \cite{octa,Goldberg} that all 
octahedrites of symmetry 
$O$ or $O_h$ are obtained from the 
Goldberg-Coxeter construction, i.e. they are of the form $GC_{k,l}(Octahedron)$
for some integer $0\leq l\leq k$.
The pairs $(k,l)$ correspond to the relative position of the triangles; see
Figure \ref{FirstExample} for the 
smallest such graphs and \cite{Goldberg} for
more details on the construction itself.

\begin{theorem}
All $4$-self-hedrites of symmetry $T$ or $T_d$ are of the form 
$Med^{-1}(GC_{k,l}(Octahedron))$ with $k+l$ odd.
\end{theorem}
\proof If $G$ is a $4$-self-hedrite of symmetry $T$ or $T_d$ then
its medial $G'=Med(G)$ is an octahedrite of symmetry $O$ or $O_h$.
So, $G'=GC_{k,l}(Octahedron)$ for some $(k,l)$.
The automorphism group
of the plane graph $GC_{k,l}(Octahedron)$ is transitive on triangles; so, we only need 
to determine when the triangles are not all in ${\mathcal F}_1(G')$ or
${\mathcal F}_2(G')$.
Clearly, this correspond to $k+l$ odd. \qed

For a plane graph $G$ a {\em zigzag} is a circuit of edges, such that any two
but no three, consecutive edges belong to the same face. 
Zigzags of $G$ correspond to central circuits of $Med(G)$, see 
an example on Figure \ref{ExampleZigzagCCmedial}.
So, if $G$ is a $4$-self-hedrite with simple 
zigzags, then $Med(G)$ is an octahedrite with simple 
central circuits.
By Section \ref{CentralCircuit}, octahedrites $G'$ with simple central circuits are obtained by taking the ones of 
Figure \ref{IrreducibleSimpleOcta} and splitting each central circuit
$C_i$ into $m_i$ parallel central circuits.
Then we have to determine for which $m=(m_i)$ the triangles are in two
parts ${\mathcal F}_1(G')$ and ${\mathcal F}_2(G')$ which are equivalent
under an automorphism of $G'$. This requires a detailed analysis of the
automorphism and a search of the necessary relations between $m_i$ and parity
conditions. The details are very cumbersome but in principle we can get a
classification of the $4$-self-hedrites with simple zigzags.

In particular, 1, 3, 4, 5, 6, 7th irreducible octahedrites in Figure \ref{IrreducibleSimpleOcta} are the
medial graphs of 1, 6, 7, 11, 13, 16th $4$-self-hedrites in Figure \ref{MinimalRepresentativeSelfDual34}, respectively; they are all {\em irreducible} 
$4$-self-hedrites with simple zigzags.

\begin{figure}
\begin{center}
\begin{minipage}[b]{2.0cm}
\centering
\resizebox{0.5cm}{!}{\rotatebox{90}{\includegraphics{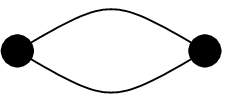}}}\par
$D_{2h}$, $n=2$
\end{minipage}
\begin{minipage}[b]{2.0cm}
\centering
\resizebox{0.35cm}{!}{\rotatebox{90}{\includegraphics{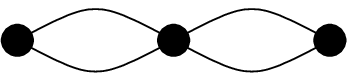}}}\par
$C_{2v}$, $n=3$
\end{minipage}
\begin{minipage}[b]{2.0cm}
\centering
\resizebox{1.0cm}{!}{\rotatebox{90}{\includegraphics{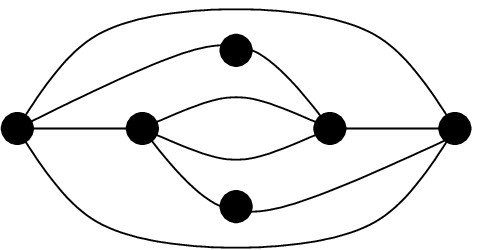}}}\par
$D_{2}$, $n=6$
\end{minipage}
\begin{minipage}[b]{2.0cm}
\centering
\resizebox{2.0cm}{!}{\rotatebox{0}{\includegraphics{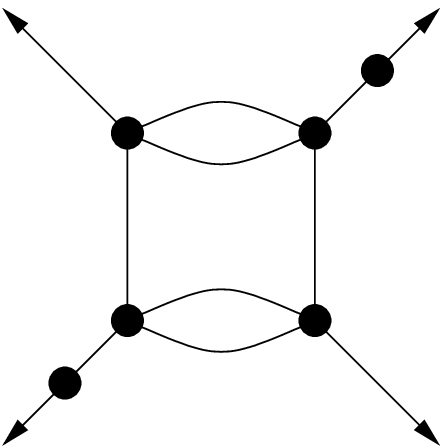}}}\par
$C_{2}$, $n=7$
\end{minipage}
\begin{minipage}[b]{3.0cm}
\centering
\resizebox{2.8cm}{!}{\rotatebox{0}{\includegraphics{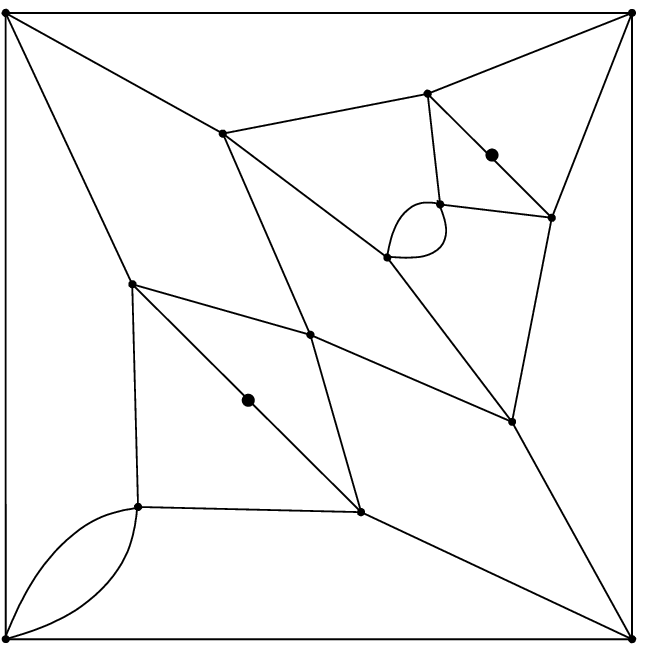}}}\par
$C_{2h}$, $n=16$
\end{minipage}
\end{center}
\caption{Minimal representatives for each possible symmetry group of $2$-self-hedrites}
\label{MinimalRepresentative2SelfHedrites}
\end{figure}

\begin{figure}
\begin{center}
\begin{minipage}[b]{3.0cm}
\centering
\resizebox{2.8cm}{!}{\rotatebox{0}{\includegraphics{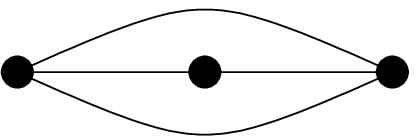}}}\par
$C_{2v}$, $n=3$
\end{minipage}
\begin{minipage}[b]{3.0cm}
\centering
\resizebox{2.8cm}{!}{\rotatebox{90}{\includegraphics{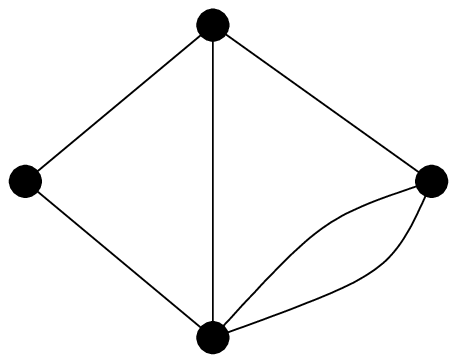}}}\par
$C_{1}$, $n=4$
\end{minipage}
\begin{minipage}[b]{3.0cm}
\centering
\resizebox{2.8cm}{!}{\rotatebox{90}{\includegraphics{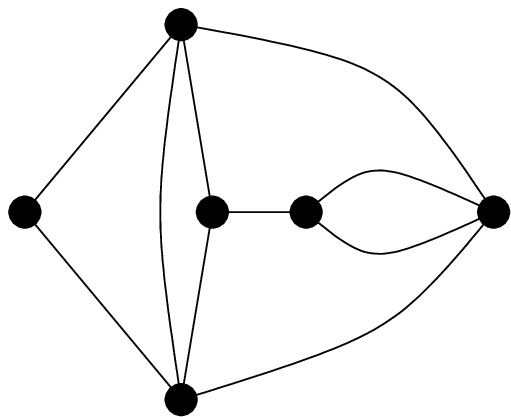}}}\par
$C_{s}$, $n=6$
\end{minipage}
\begin{minipage}[b]{3.0cm}
\centering
\resizebox{2.8cm}{!}{\rotatebox{0}{\includegraphics{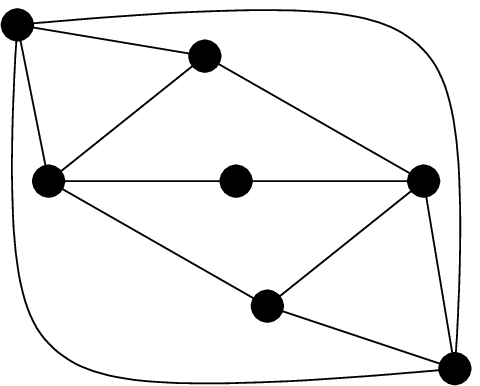}}}\par
$C_{2}$, $n=7$
\end{minipage}
\end{center}
\caption{Minimal representatives for each possible symmetry group of $3$-self-hedrites}
\label{MinimalRepresentative3SelfHedrites}
\end{figure}

\begin{figure}
\begin{center}
\begin{minipage}[b]{2.0cm}
\centering
\epsfig{height=16mm, file=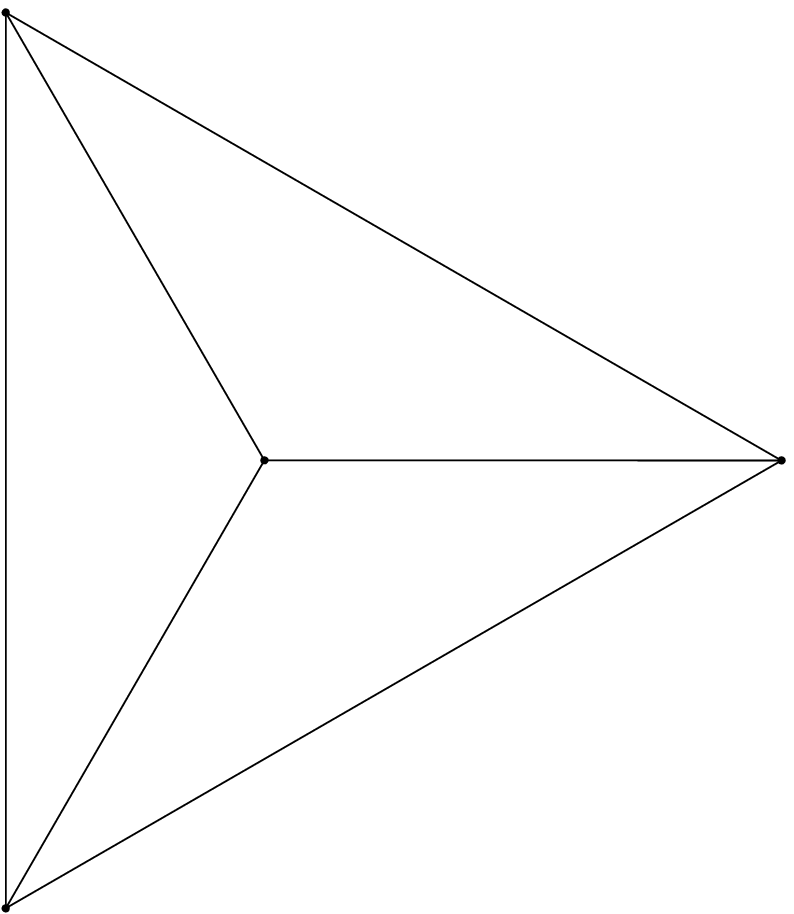}\par
$T_{d}$, $n=4$
\end{minipage}
\begin{minipage}[b]{2.3cm}
\centering
\epsfig{height=16mm, file=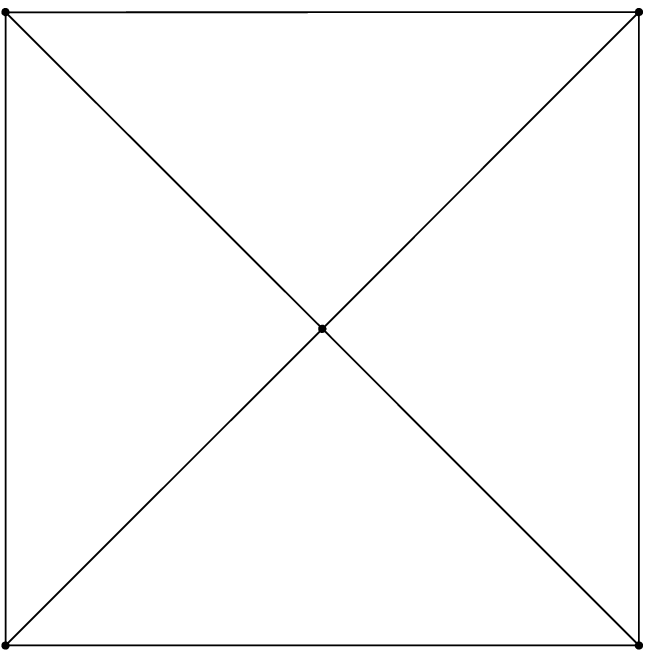}\par
$C_{4\nu}$, $n=5$
\end{minipage}
\begin{minipage}[b]{2.3cm}
\centering
\epsfig{height=16mm, file=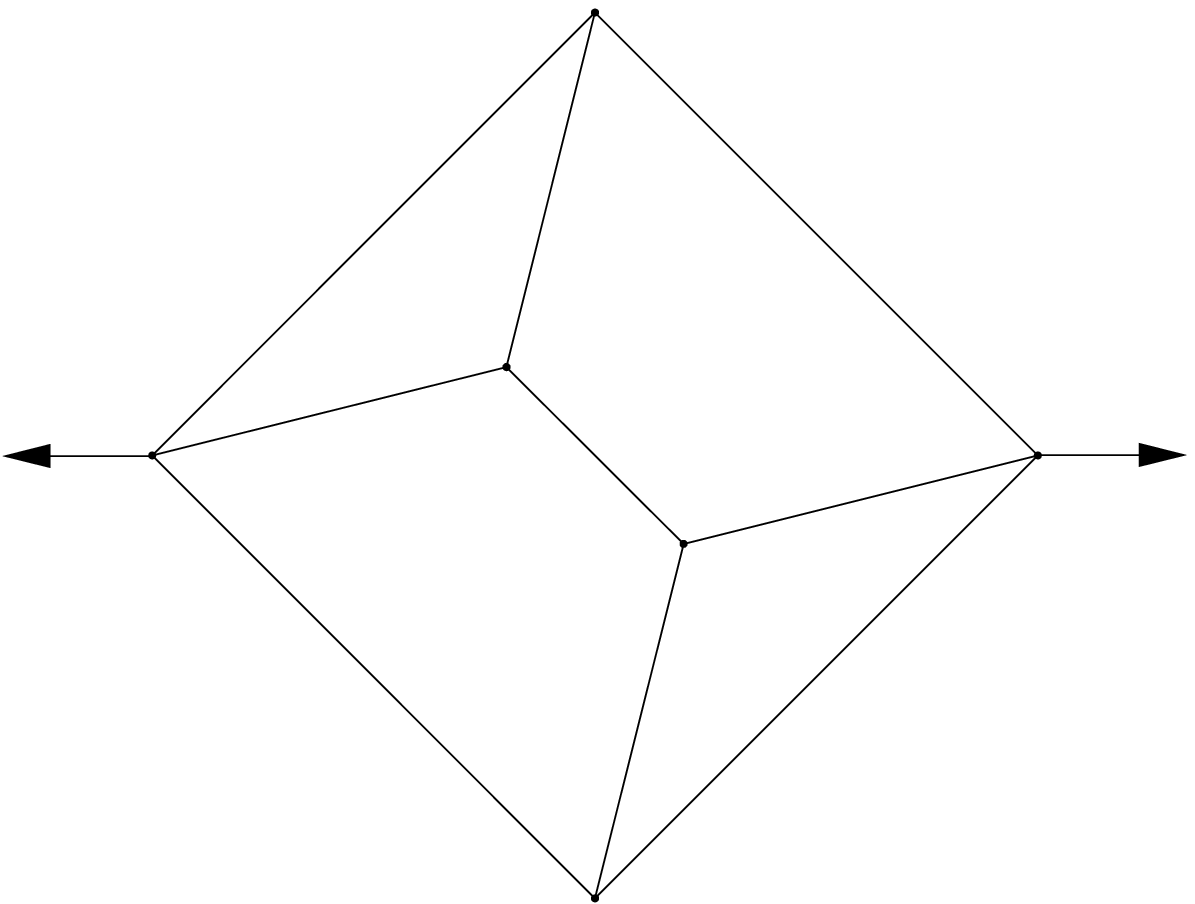}\par
$C_2$, $n=6$
\end{minipage}
\begin{minipage}[b]{2.0cm}
\centering
\epsfig{height=16mm, file=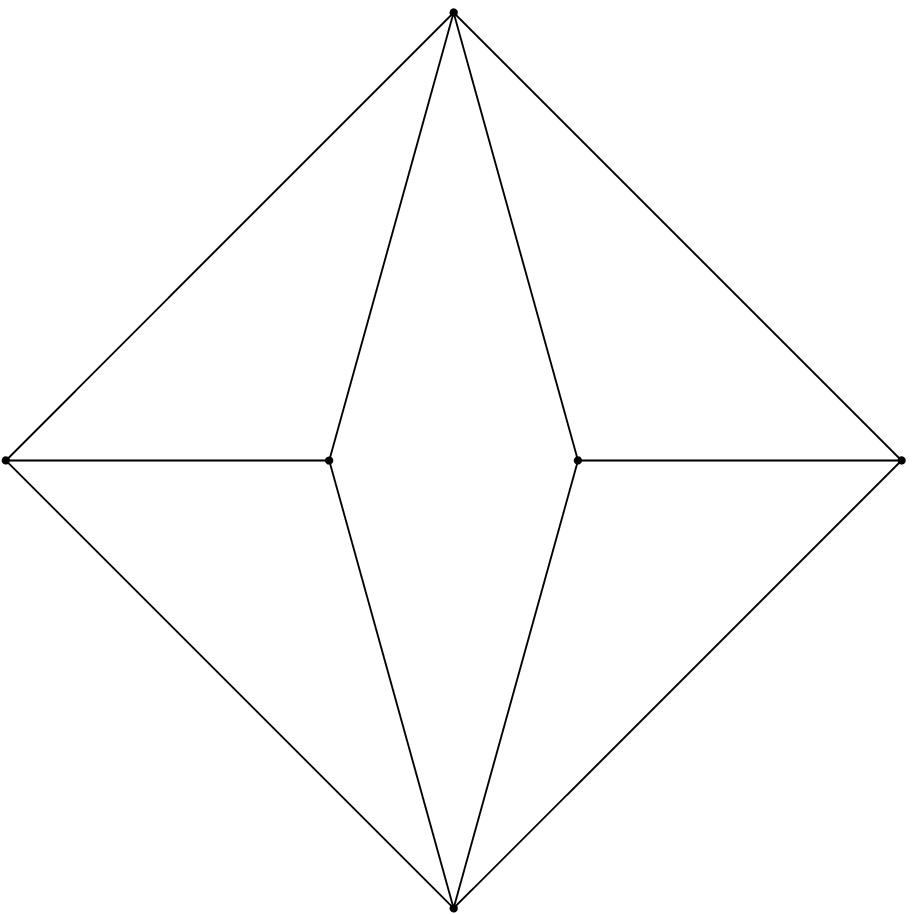}\par
$D_{2h}$, $n=6$
\end{minipage}
\begin{minipage}[b]{2.3cm}
\centering
\epsfig{height=16mm, file=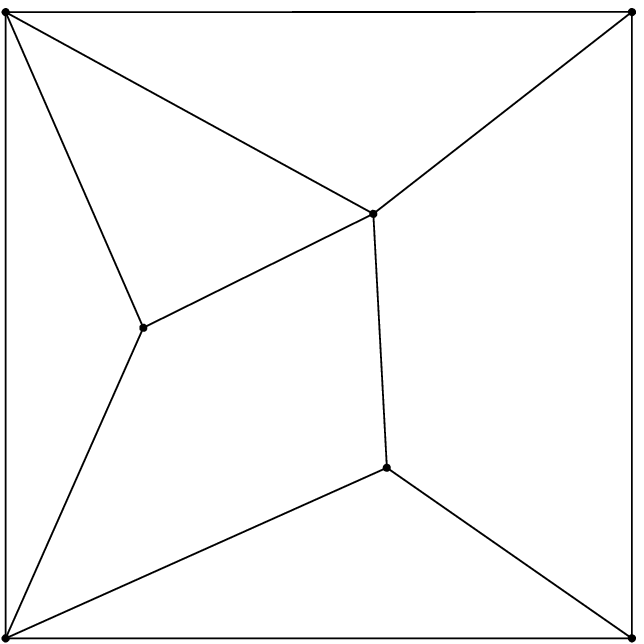}\par
$C_1$, $n=7$
\end{minipage}
\begin{minipage}[b]{2.3cm}
\centering
\epsfig{height=16mm, file=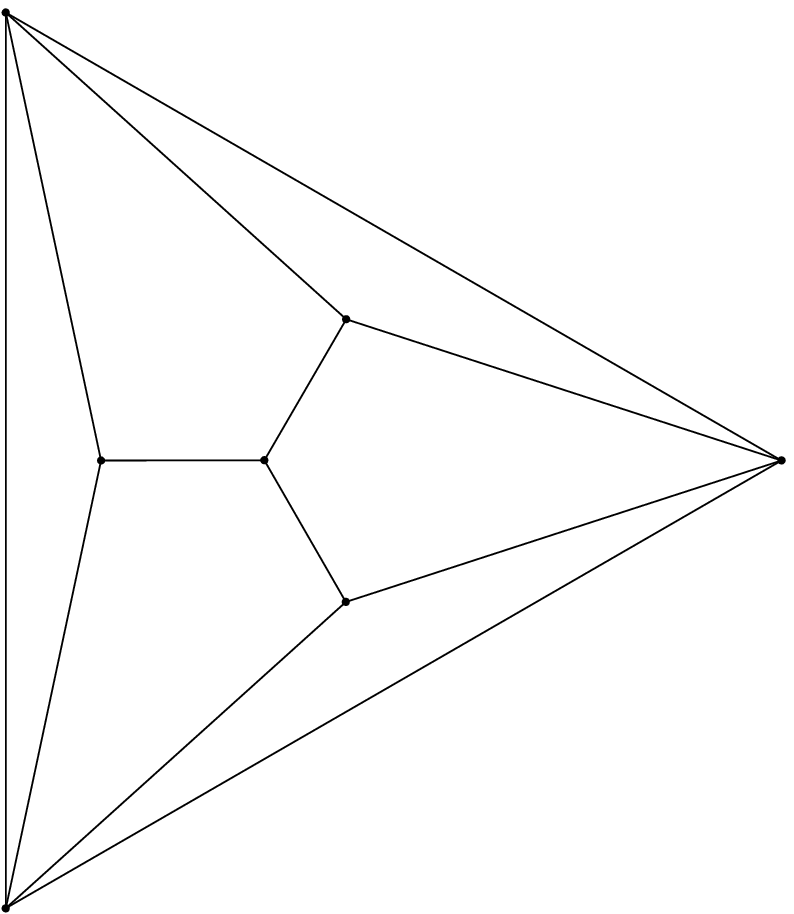}\par
$C_{3\nu}$, $n=7$
\end{minipage}
\begin{minipage}[b]{2.0cm}
\centering
\epsfig{height=16mm, file=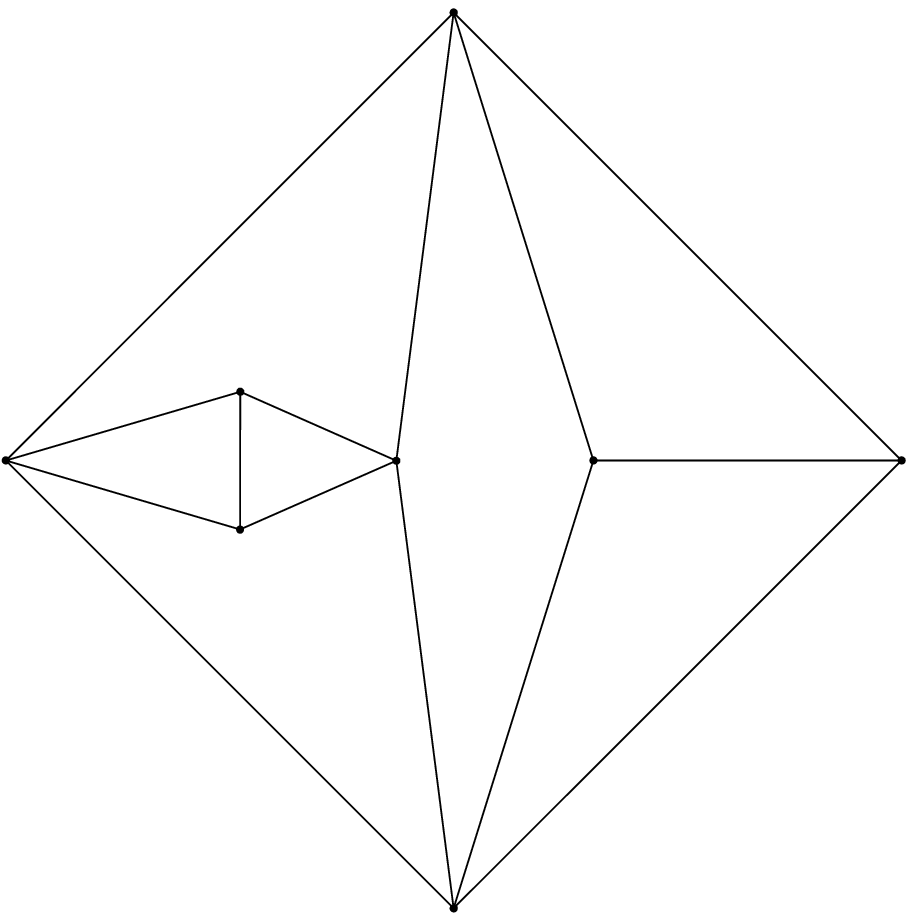}\par
$D_{2d}$, $n=8$
\end{minipage}
\begin{minipage}[b]{2.3cm}
\centering
\epsfig{height=16mm, file=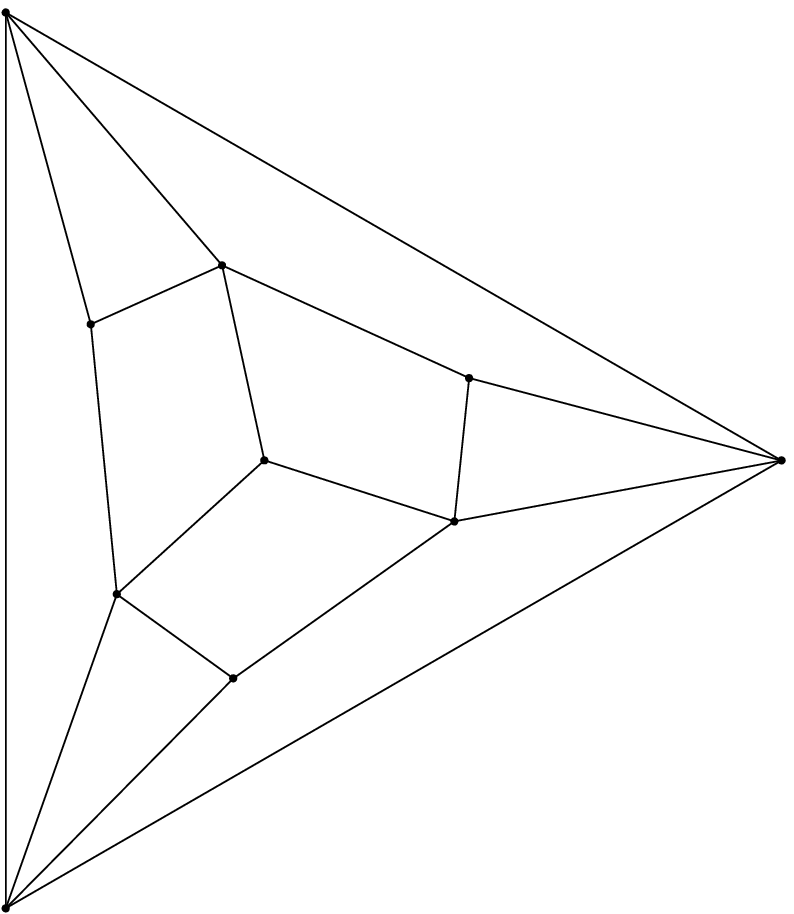}\par
$C_3$, $n=10$
\end{minipage}
\begin{minipage}[b]{2.3cm}
\centering
\epsfig{height=16mm, file=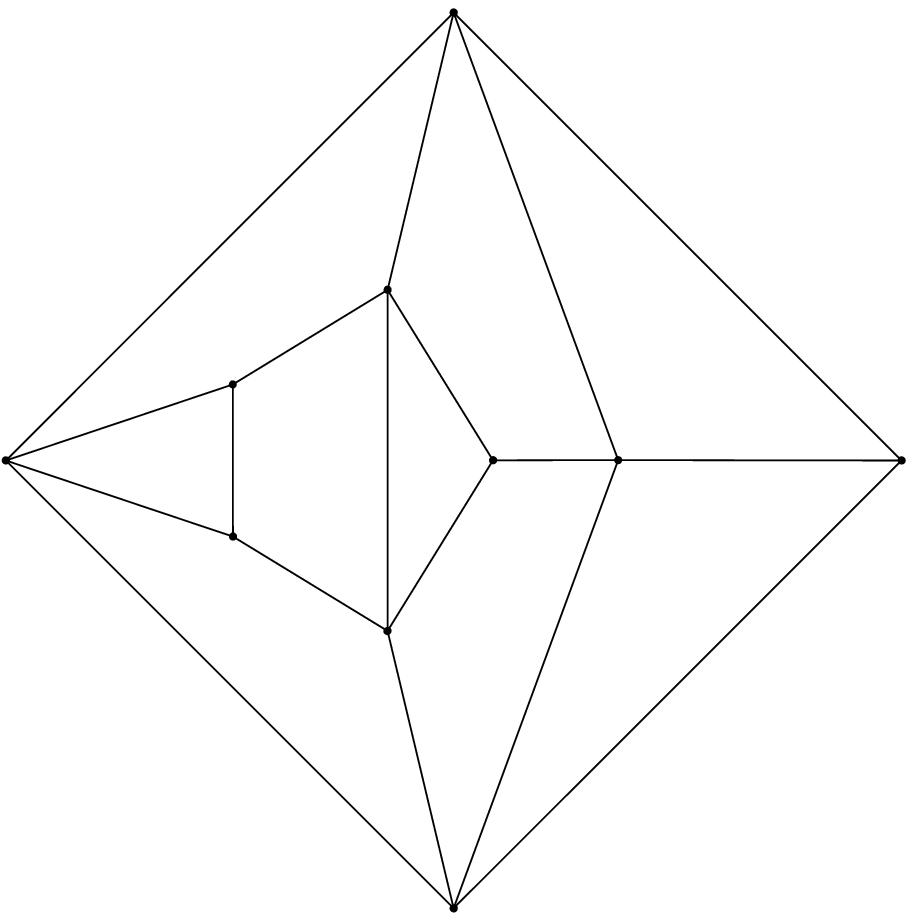}\par
$C_{s}$, $n=10$
\end{minipage}
\begin{minipage}[b]{2.0cm}
\centering
\epsfig{height=16mm, file=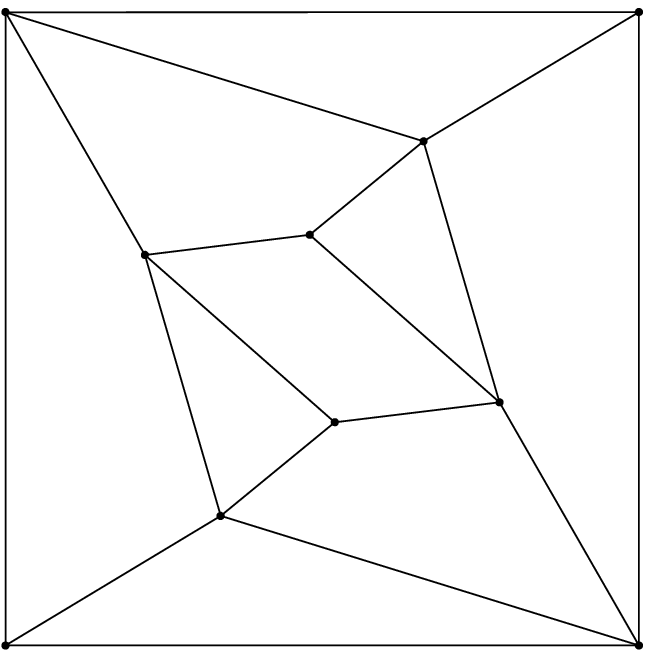}\par
$D_{2}$, $n=10$
\end{minipage}
\begin{minipage}[b]{2.3cm}
\centering
\epsfig{height=16mm, file=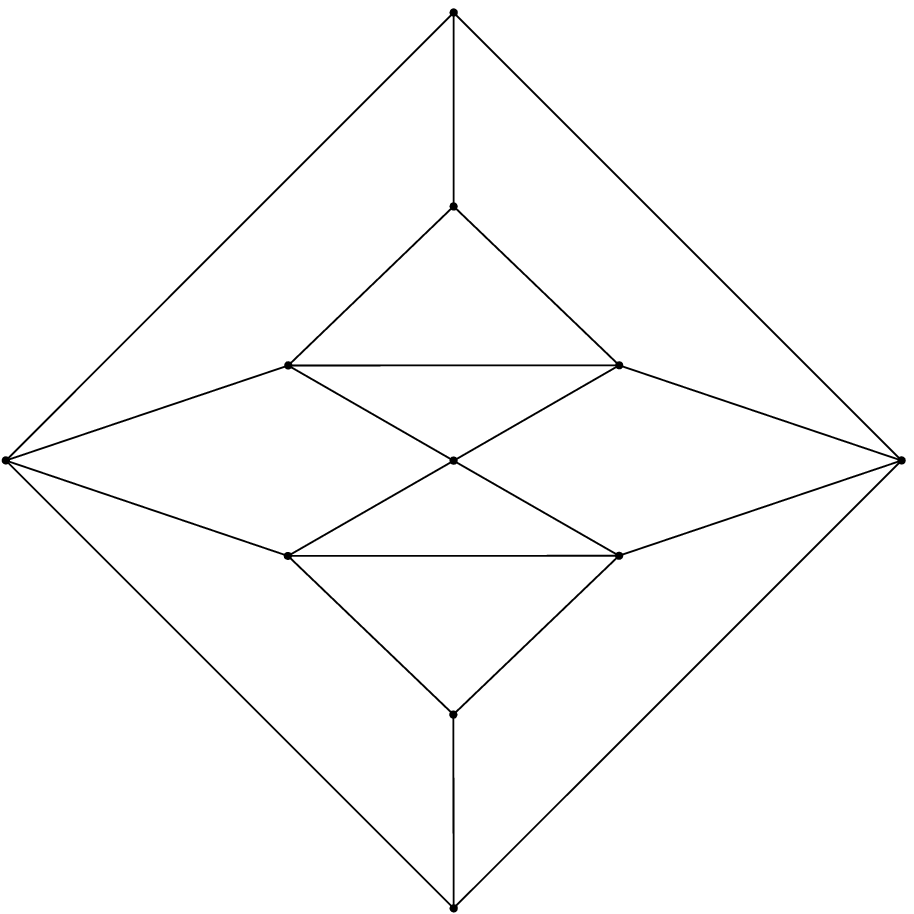}\par
$C_{2\nu}$, $n=11$
\end{minipage}
\begin{minipage}[b]{2.3cm}
\centering
\epsfig{height=16mm, file=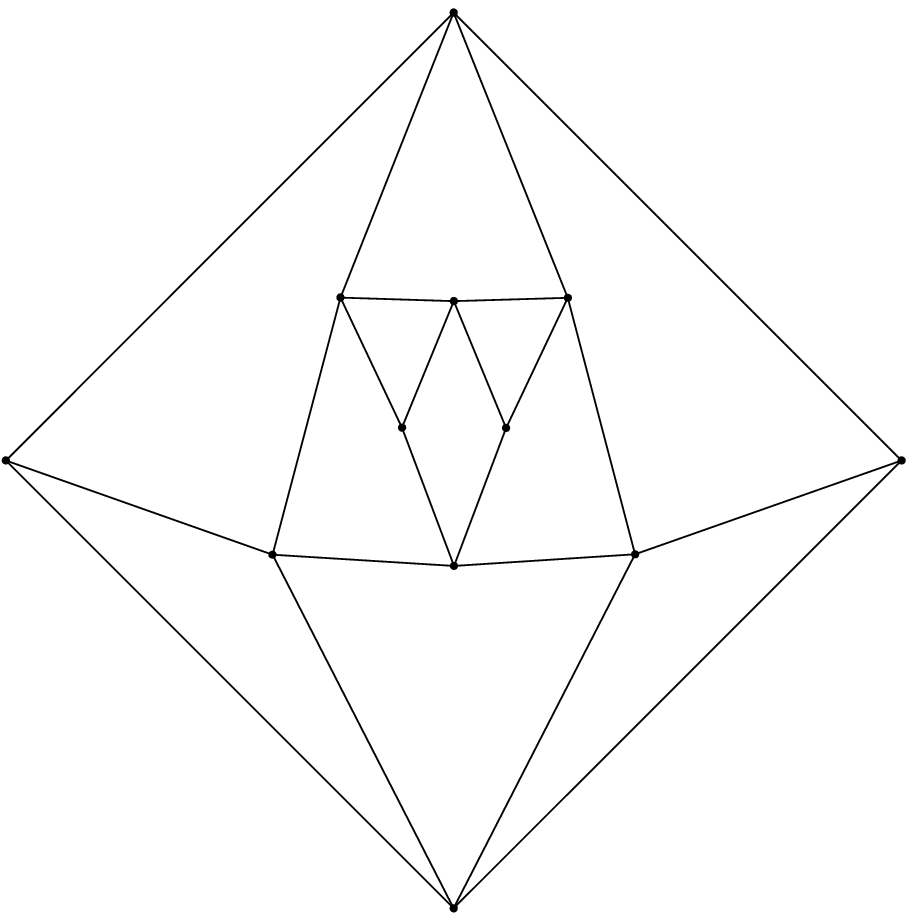}\par
$C_{2h}$, $n=12$
\end{minipage}
\begin{minipage}[b]{2.0cm}
\centering
\epsfig{height=16mm, file=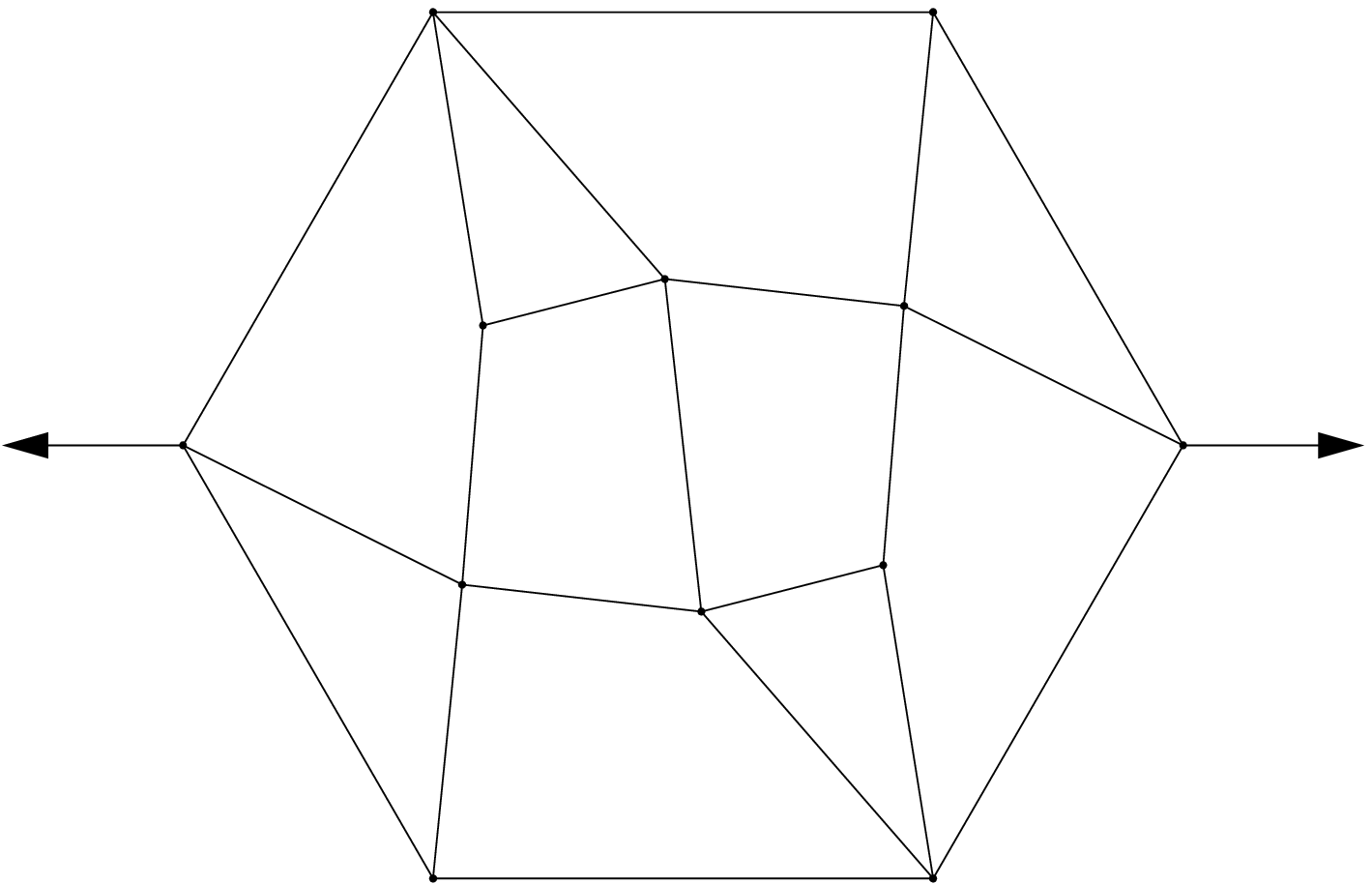}\par
$S_{4}$, $n=12$
\end{minipage}
\begin{minipage}[b]{2.3cm}
\centering
\epsfig{height=16mm, file=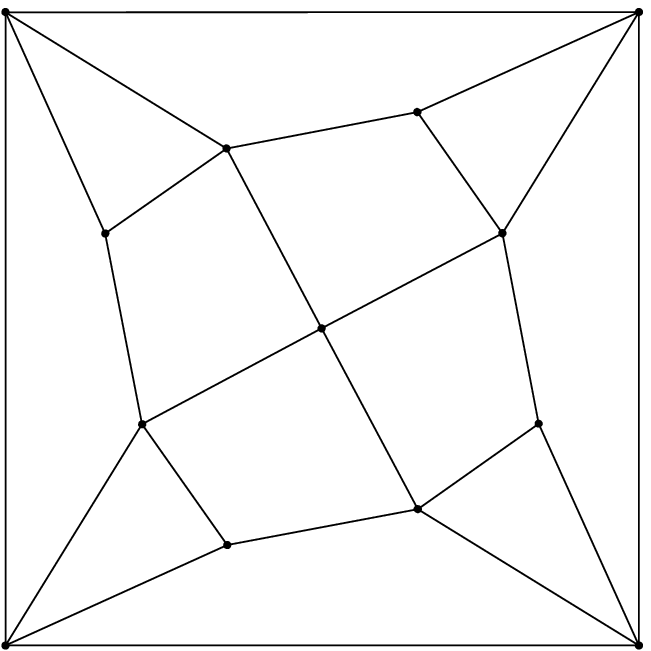}\par
$C_4$, $n=13$
\end{minipage}
\begin{minipage}[b]{2.3cm}
\centering
\epsfig{height=16mm, file=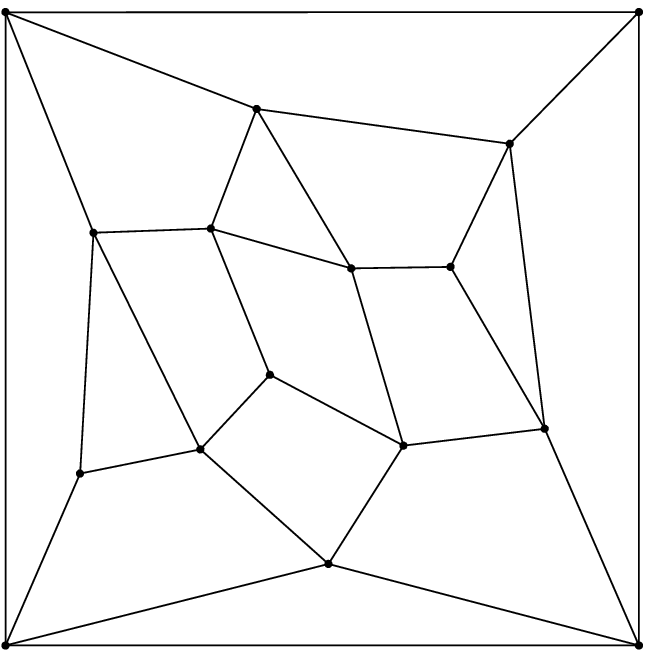}\par
$C_{i}$, $n=16$
\end{minipage}
\begin{minipage}[b]{2.0cm}
\centering
\epsfig{height=16mm, file=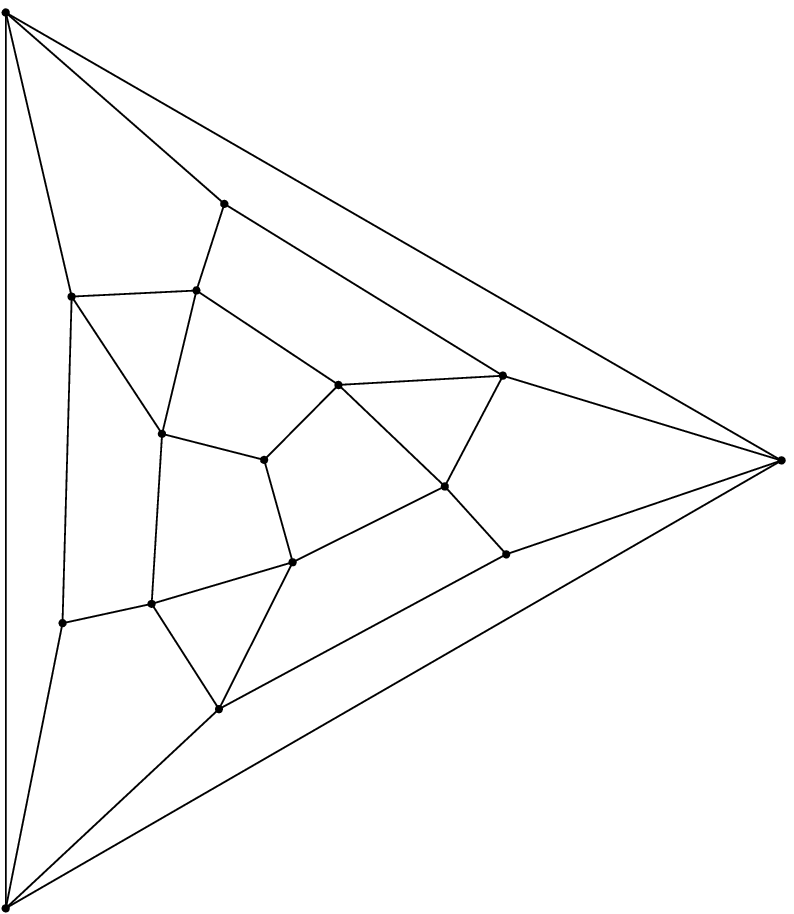}\par
$T$, $n=16$
\end{minipage}

\end{center}
\caption{Minimal representatives for each possible symmetry group of $4$-self-hedrites}
\label{MinimalRepresentativeSelfDual34}
\end{figure}

\begin{figure}

\begin{center}
\begin{minipage}{4.0cm}
\centering
\epsfig{height=24mm, file=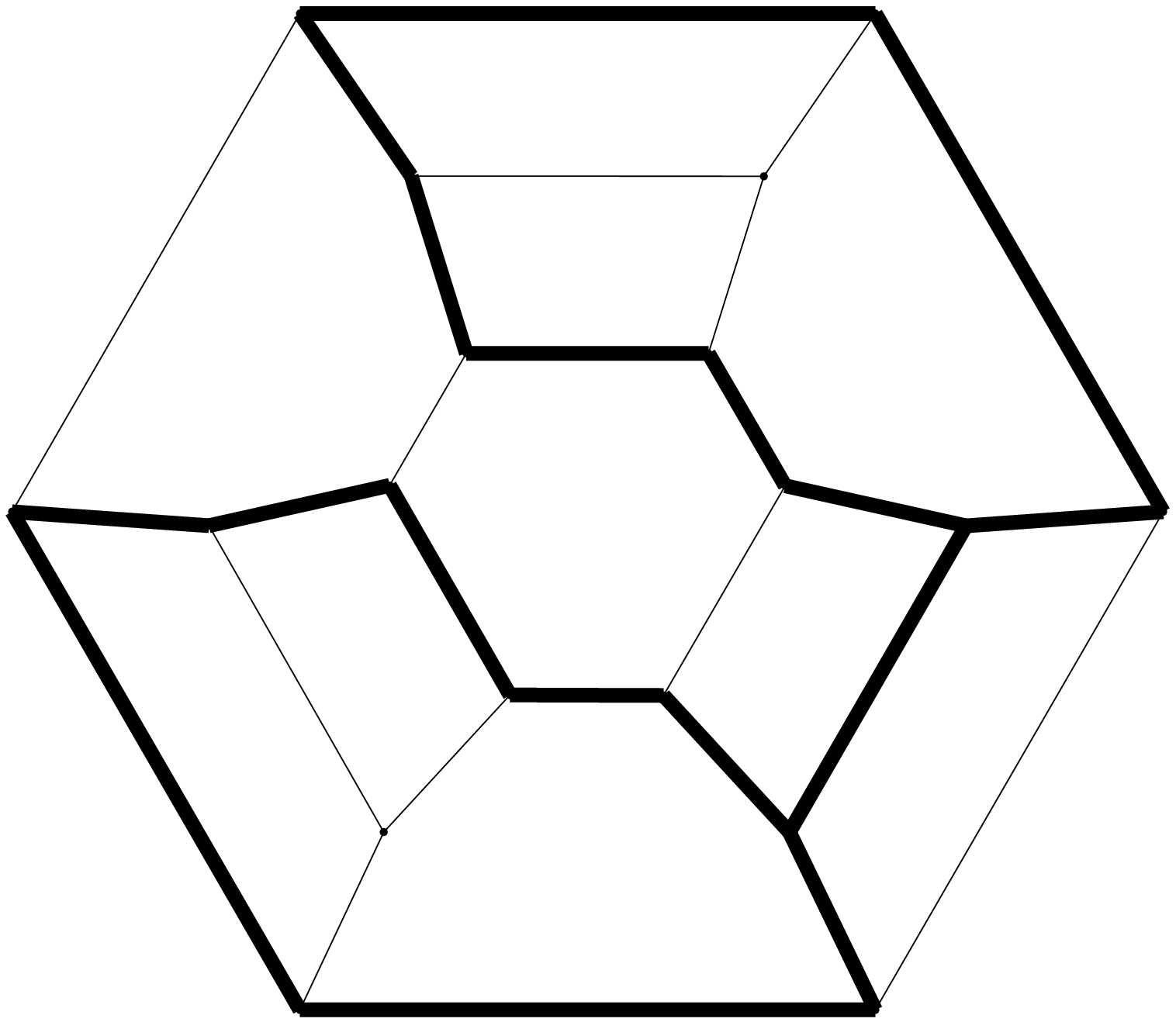}\par
$G$
\end{minipage}
\begin{minipage}{5.6cm}
\centering
\epsfig{height=24mm, file=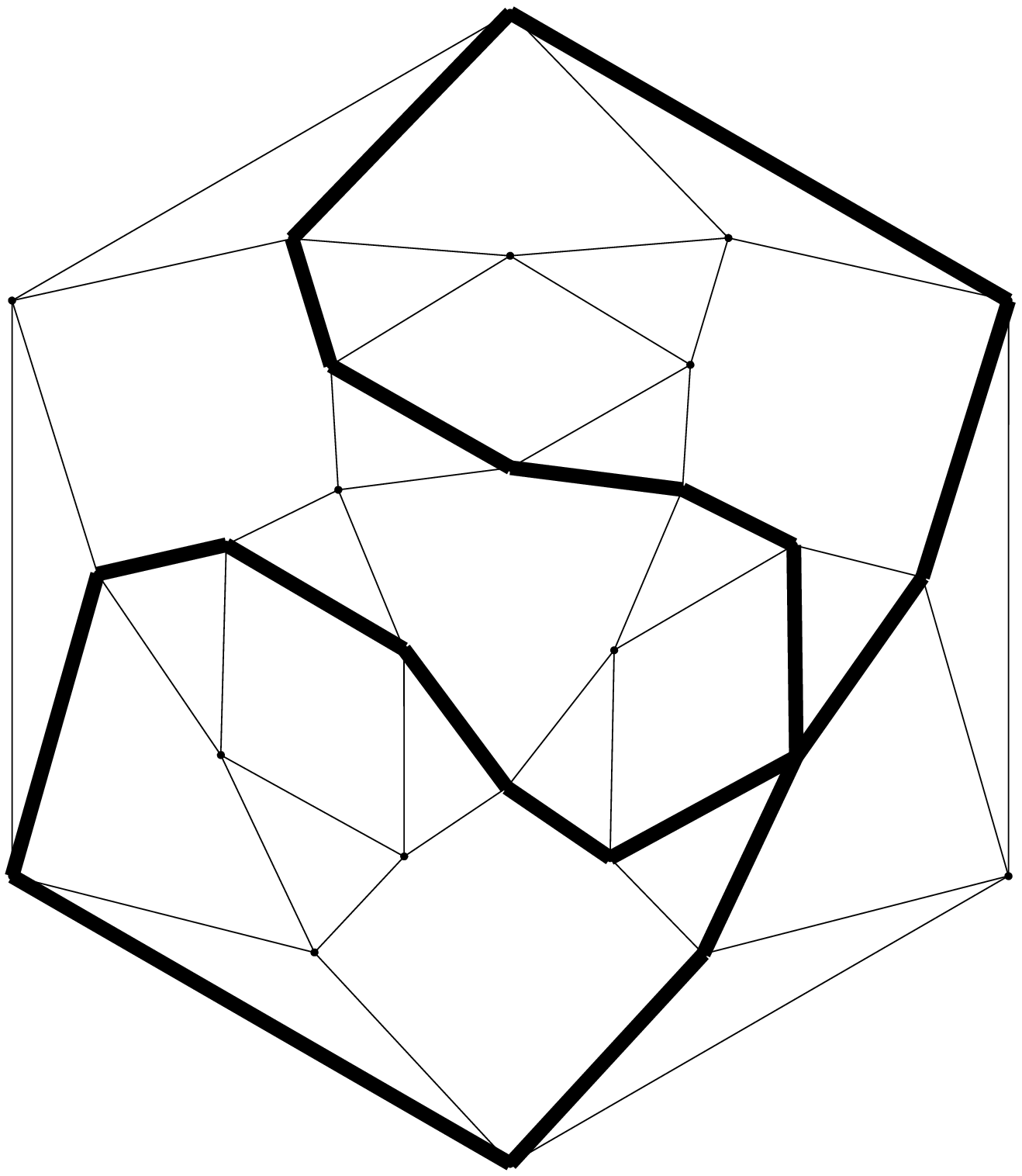}\par
$Med(G)$
\end{minipage}
\end{center}
\caption{Example of a zigzag in a plane graph $G$ and the corresponding central circuit in $Med(G)$}
\label{ExampleZigzagCCmedial}
\end{figure}

\section{Going on surfaces}
In \cite{surfaces} was considered a generalization of plane fullerenes on any irreducible 
surface.
Similarly, it is easy to see that any {\em generalized octahedrite}, i.e., a $4$-regular map on 
 an irreducible surface, having only $3$- and $4$-gonal faces, is either an octahedrite on 
sphere $S^2$, or a partition of torus $T^2$ (or
Klein bottle $K^2$) into $4$-gons, or the antipodal quotient of a centrally
symmetric octahedrite on the projective plane $P^2$ (having $4$ $3$-gonal
faces).
Maps on surfaces of high genus can be very complicated. Actually, there
are examples with the genus being about half the number of vertices.
Here, for the sake of simplicity and the search of more complex examples,
we limit ourselves to graphs with no loops or multiple edges.
The {\em minimal}, i.e. with minimal number of vertices, generalized
octahedrite on $S^2$ is Octahedron $K_{2,2,2}$; on $P^2$ it is the antipodal
quotient of Cube with two opposite faces triangulated in their center, that is
$K_5$.
On $T^2$ it is $K_5$, and on $K^2$ it is 
again $K_{2,2,2}$ (but embedded as a quadrangulation); see Fig. 6 in \cite{Nakamoto}.

Finally, it is easy to 
check that any {\em generalized $4$-self-hedrite}, i.e., self-dual 
 map on an irreducible surface, having only $3$- and $4$-gonal faces,
is either a $4$-self-hedrite on sphere $S^2$, or a $4$-regular partition of torus $T^2$ (or
Klein bottle $K^2$) into $4$-gons, or the antipodal quotient of a centrally
symmetric $4$-self-hedrite on the projective plane $P^2$ (having $2$ $3$-gonal faces).
The minimal generalized 
$4$-self-hedrite graph on $S^2$ is Tetrahedron; on $P^2$ it is the
antipodal quotient of the 12th graph on Figure \ref{MinimalRepresentativeSelfDual34}, that is complete graph $K_6$ with disjoint $2$- and $4$-vertex paths deleted.
On $T^2$ it is $K_5$, and on $K^2$ it is again $K_{2,2,2}$
(see Fig. 6 in \cite{Nakamoto}) embedded as a quadrangulation.

Similar results hold for generalization of $i$-hedrites and $i$-self-hedrites
from sphere on any irreducible surface.
On $T^2$ and $K^2$ it gives the $4$-regular quadrangulations.
On $P^2$ they are the antipodal quotients of such centrally symmetric graphs on $S^2$.
So, $2p_2+p_3$ becomes $4$ for $i$-hedrites and $2$ for $i$-self-hedrites on $P^2$.

\section{Acknowledgment}
First author has been supported by the Croatian Ministry of Science, Education and Sport under contract 098-0982705-2707.
The authors thank G. Brinkmann for help with the {\tt ENU} program.

\end{document}